\documentclass[11pt]{amsart}
\usepackage{lmodern}
\usepackage[T1]{fontenc}
\usepackage{microtype}
\usepackage{amssymb}
\usepackage{amsthm}
\usepackage{amscd}
\usepackage{hyperref}
\usepackage{cite}
\usepackage{verbatim}

\newtheorem{theorem}{Theorem}[section]
\newtheorem{lemma}[theorem]{Lemma}

\newtheorem{cor}[theorem]{Corollary}

\newtheorem*{theorem*}{Theorem}
\newtheorem*{cor*}{Theorem}

\theoremstyle{definition}
\newtheorem{definition}[theorem]{Definition}
\newtheorem{notation}[theorem]{Notation}
\newtheorem{remark}[theorem]{Remark}

\newcounter{claimcounter}
\numberwithin{claimcounter}{theorem}

\makeatletter

\def\dotminussym#1#2{%
  \setbox0=\hbox{$\m@th#1-$}%
  \kern.5\wd0%
  \hbox to 0pt{\hss\hbox{$\m@th#1-$}\hss}%
  \raise.6\ht0\hbox to 0pt{\hss$\m@th#1.$\hss}%
  \kern.5\wd0}

\mathchardef\mhyphen="2D

\allowdisplaybreaks[2]

\newcommand{\GB}{\text{Gr\"{o}bner }}

\begin{document}


\title[Polynomial bounds in polynomial rings]{Explicit polynomial bounds on prime ideals in polynomial rings over fields}
\author{William Simmons and Henry Towsner}
\date{\today}
\thanks{Partially supported by NSF grant DMS-1600263.}
\thanks{The authors thank Florian Pop, Thomas Preu, and Matthias Aschenbrenner for helpful discussions.}
\address{Department of Mathematics and Computer Science, Hobart and William Smith Colleges, 300 Pulteney Street, Geneva, NY 14456, USA}
\email{wsimmons@hws.edu}
\urladdr{\url{https://math.hws.edu/simmons/}}
\address {Department of Mathematics, University of Pennsylvania, 209 South 33rd Street, Philadelphia, PA 19104-6395, USA}
\email{htowsner@math.upenn.edu}
\urladdr{\url{http://www.math.upenn.edu/~htowsner}}

\begin{abstract}
Suppose $I$ is an ideal of a polynomial ring over a field, $I\subseteq k[x_1,\ldots,x_n]$, and whenever $fg\in I$ for $f,g$ of degree $\leq b$, then either $f\in I$ or $g\in I$.  When $b$ is sufficiently large, it turns out that $I$ is prime. Schmidt-G\"ottsch proved \cite{schmidt89} that ``sufficiently large'' can be taken to be a polynomial in the degree of generators of $I$ (with the degree of this polynomial depending on $n$).  However, Schmidt-G\"ottsch used model-theoretic methods to show this and did not give any indication of how large the degree of this polynomial is.  In this paper we obtain an explicit bound on $b$, polynomial in the degree of the generators of $I$. We also give a similar bound for detecting maximal ideals in $k[x_1,\ldots,x_n]$.
\end{abstract}

\maketitle


\section{Introduction}




Suppose $I$ is an ideal of a polynomial ring over a field, $I\subseteq k[x_1,\ldots,x_n]$, and whenever $fg\in I$ for $f,g$ of degree $\leq b$, then either $f\in I$ or $g\in I$. When $b$ is sufficiently large, it turns out that $I$ is prime.   Schmidt-G\"ottsch proved \cite{schmidt89} that ``sufficiently large'' can be taken to be a polynomial in the degree of generators of $I$ (with the degree of this polynomial depending on $n$).  However Schmidt-G\"ottsch used model-theoretic methods to show this, and did not give any indication of how large the degree of this polynomial is.  In this paper we obtain an explicit bound on $b$, polynomial in the degree of the generators of $I$. We also give a similar bound for detecting maximal ideals in $k[x_1,\ldots,x_n]$.

Our problem belongs to a thread of algorithmic algebra going back over a century. Drawing on the work of Kronecker, Noether, and others, Hermann published in 1926 a seminal paper entitled ``Die Frage der endlich vielen Schritte in der Theorie der Polynomideale'' (The question of finitely many steps in polynomial ideal theory) \cite{MR1512302} (see \cite{hermann1998question} for an English translation with commentary). Hermann obtained bounds on calculations in polynomial rings over fields such as witnessing ideal membership and finding bases of syzygies. Other tasks she showed to be algorithmically possible include primary decomposition, intersection, and quotients of ideals. 


%
In \cite{seidenberg1974constructions} Seidenberg analyzed, extended, and in some cases corrected Hermann's work. He noted (p. 311, \cite{seidenberg1974constructions}) that one can decide primality of polynomial ideals over ``explicitly given fields'' by determining if the ideal is primary and, if so, whether the ideal is equal to its associated prime. Seidenberg gave bounds on some steps of the process but did not analyze the complexity of deciding primality.

Buchberger's method of \GB bases (introduced in his thesis \cite{buchberger1965algorithmus}; see Section \ref{prelim} for basic results) provided strong impetus to the development of computational algebra. While many problems are in principle solvable by the methods of Hermann and Seidenberg, the theory of \GB bases has grown substantially and is often more convenient. In particular, much work has been done on the complexity of \GB bases. Mayr and Mayer \cite{mayr1982complexity} established doubly exponential worst-case bounds on the complexity of \GB bases (as a function of the number of indeterminates). Building on subsequent efforts, Dub\'{e} \cite{dube1990structure} gave a sharper (but still doubly exponential) bound that we take as a starting point for our analysis in \ref{thm:GB Dube}. 

\GB basis methods can be used as building blocks for algorithms that test primality in polynomial rings. Gianni, Trager, and Zacharias \cite{gianni1988grobner} (pp.~157-158) described an algorithm for deciding primality by, among other things, inductively reducing the number of variables to the univariate case. Eisenbud, Huneke, and Vasconcelos  \cite{eisenbud1992direct} gave a method for carrying out primary decomposition of ideals in polynomial rings over perfect fields; one may then follow Seidenberg's observation to decide primality. To our knowledge no one has explicitly analyzed the complexity of these algorithms. However, many of our bounds in Section \ref{prelim} concern steps of the algorithm from \cite{gianni1988grobner} as presented in \cite{adams1994introduction}.

Another result related to our problem is Lemma 10 of \cite{chistov2009double}. There Chistov shows that there is a constant $c\in \mathbb{N}$ such that, given an algebraic variety defined by polynomials in $n$ variables of total degree at most $d$, irreducible components of the variety are defined by polynomials of total degree at most $d^{2^{cn}}$.













Our main result is:
\begin{theorem}\label{thm:prime bound}
There exists a function $b(n,d)$ that is polynomially-bounded in $d$ and such that for any field $k$ and proper ideal $I\subseteq k[x_1,\ldots,x_n]$ with generators of total degree at most $d$, if $fg\in I$ implies that either $f\in I$ or $g\in I$ for all $f,g$ of degree less than or equal to $b(n,d)$, then $I$ is prime. Moreover, $b(n,d)\leq max\left\{ (2d)^{2^{3n^2+2n}}, (2nd)^{n(n+1)^2{2^{4n+6}}}\right\}$ for all $n,d\geq 1$. If $n=1,d>1$ or $d=1$ the values $b(1,d)=d$ or $b(n,1)= 0$, respectively, suffice to ensure primality.

\end{theorem}

\noindent See Theorems \ref{thm:prime} and \ref{thm:final estimate} for the final steps of the proof. The function $b(n,d)$ is defined in the course of the argument and a more compact bound easily follows (see \ref{thm:simple bd}):
\begin{cor}
$b(n,d)\leq (2nd)^{n^32^{6n^2}}$ for all $n,d\geq 1$.
\end{cor}

Schmidt-G\"ottsch \cite{schmidt89} showed that the function $b(n,d)$ is polynomially-bounded in $d$ using ultraproduct methods (as in \cite{vdDS}, \cite{schoutens2010use}, and \cite{Goral18}), but did not give any explicit calculation of the bound.  Abstract results from proof theory \cite{kohlenbach:MR2445721,1804.10809} suggest that it is possible to extract explicit bounds from proofs that use ultraproducts.  The authors demonstrated \cite{simmons_towsner} that such methods can be used to extract explicit bounds from ultraproduct proofs in commutative and differential algebra \cite{HTKM}\cite{vdDS}.  (For instance, in Theorem 2.8 and Lemma 7.7 of \cite{simmons_towsner} we gave another bound on $b(n,d)$, but that bound was not polynomial in the degree of the generators.)  However the arguments here, while inspired by Schmidt-G\"ottsch's approach, were not purely extracted by proof theory; instead we bring in known results using Gr\"obner bases, etc., where possible and add new arguments to simplify (or clarify) the most complicated parts of Schmidt-G\"ottsch's argument.



\begin{section}{Outline of the argument}
We briefly describe our strategy. We are interested in \emph{counterexamples to primality}; that is, polynomials $f,g\notin I$ such that $fg\in I$. Our goal is to show that if there are no counterexamples to primality of $I$ having degree  $b(n,d)$ 
or less (we say $I$ is \emph{prime up to} $b(n,d)$), then there are no counterexamples of any degree (i.e., $I$ is prime). The main ingredients are  localization of polynomial rings, bounds on Gr\"{o}bner bases for various auxiliary ideals, and bounds on solutions to systems of linear equations over polynomial rings. 

\begin{itemize}
\item Relabeling if necessary, choose a maximal set of indeterminates $x_1,\ldots, x_r$ that is algebraically independent modulo $I$. By independence mod $I$ we have $k[x_1,\ldots,x_r]\cap I =\{0\}$; by maximality $k[x_1,\ldots, x_r, x_j]\cap I\neq \{0\}$ for any $j>r$.    Without loss of generality, we may assume that no $x_i\in I$ (otherwise we solve the problem for smaller $n$).

\item Localize with respect to the variables $x_1,\ldots,x_r$; that is, consider the ideal $J=Ik(x_1,\ldots,x_r)[x_{r+1},\ldots, x_n]$ generated by $I$ in the ring $k(x_1,\ldots,x_r)[x_{r+1},\ldots,x_n]$. The first main task is to show that $k[x_1,\ldots, x_n]\cap J=I$.  This step requires induction to prove that there are no counterexamples to primality of $I$ in which one factor belongs to $k[x_1,\ldots,x_r]$. Here we use bounds on \GB bases of saturation and quotient ideals.
\item  The second major task is to show that $J$ is maximal. We introduce a maximal ideal $M$ extending $J$ and prove that $J=M$. This step uses various tools such as the primitive element theorem, flat extensions, and bounds on \GB bases for elimination ideals. We also require some complicated inductive arguments on the number of variables $x_{r+1},\ldots, x_n$. Once we know that $J$ is maximal and hence prime, we conclude that $I$ is prime.


\end{itemize}

We are principally concerned with the inductive definition of $b(n,d)$ and obtaining a relatively clean estimate rather than precise values of constants and lower-order terms. For this reason (and for readability) our estimate of the bound is somewhat relaxed compared to the actual content of the proofs. By following closely, the interested reader may keep track of such details.
\end{section}

\section{Preliminary definitions and results}\label{prelim}
 
Throughout the section, $k$ denotes an arbitrary field (except for Theorems \ref{thm:flat bd pth power} and \ref{thm:baby fflat bd pth power}, which require that $k$ have positive characteristic).





We use the following algebraic result to control degrees when going from rational functions to polynomials:

\begin{theorem} (Gauss' Lemma)\label{thm:Gauss lemma}
Let $R$ be a Unique Factorization Domain (UFD) and let $F$ be the fraction field of $R$. If $p(x)\in R[x]$ factors as $q_1(x)q_2(x)$ in $F[x]$ such that $q_1,q_2$ have degrees $d_1,d_2$ in $x$, then $p(x)$ factors as $p_1(x)p_2(x)$ in $R[x]$ where the degrees of $p_1,p_2$ are also $d_1,d_2$, respectively.
\end{theorem}

\begin{proof}
See, for example, Corollary 2.2 in \cite{lang2002algebra}.
\end{proof}

The following version of a well-known theorem in field theory plays a prominent role in Section \ref{sec:main theorem} when we consider the case of zero-dimensional ideals. 

\begin{theorem}{(Primitive element theorem)} \label{thm:Primitive element}
Let $F\subseteq K$ be a field extension, with $F$ being an infinite field. Let $a,b \in K$ be algebraic of degrees $r,s$, respectively, over $F$. Further assume that $b$ is separable over $F$. Then for all but at most $r(s-1)$ choices of $c\in F$, the linear combination $a+cb$ generates the field $F(a,b)$ over $F$.
\end{theorem}

\begin{proof}
See \cite{van1970algebra}, p. 139. Though the bound is clear from the proof, we summarize the argument here because it is reappears in the proof of Theorem \ref{thm:separable max induction}. Let $\alpha_1,\ldots,\alpha_r$ and $\beta_1,\ldots, \beta_s$ be the roots (counting multiplicity) of the minimal polynomials $f$ of $a$ and $g$ of $b$, respectively, over $F$. (We work in an extension of $K$ containing the roots.) We may assume that $\alpha_1=a$ and $\beta_1=b$. By separability, $\beta_i\neq \beta_j$ for $i\neq j$; the list $\alpha_1,\ldots, \alpha_r$ may contain repetition. Consider the $r(s-1)$-many equations $\alpha_i + x\beta_j=a+ xb$, where $1\leq i\leq r$ and $j\neq 1$. Since $\beta_j\neq b$, each equation has at most one solution in $F$ and, except for at most $r(s-1)$ elements, none of the infinitely many remaining elements $c\in F$ satisfies any of the equations.

Let $c\in F$ satisfy none of the equations; we claim that $F(a,b)=F(a+cb)$. Consider the polynomials $g(x)$ and $f((a+cb) -cx)$, both of which belong to $F(a+cb)[x]$ and have a common root $b$. They have no other root in common: another root $\beta_j$ of $g$, with $j\neq 1$, yields an argument $(a+cb) -c\beta_j$ of $f$ that cannot equal $\alpha_i$ for any $i$ lest $c$ satisfy a forbidden equation. It follows that $g$ and $f((a+cb) -cx)$ have a GCD of degree 1 in $F(a+cb)[x]$ (if they were coprime, they could have no common root; if the GCD had degree greater than 1, they would have another common root or $g$ would have a repeated root, contradicting separability). This GCD must be $x-b$, whence $b\in F(a+cb)$ and also $a\in F(a+cb)$.
\end{proof}

\subsection{Bounds on solutions to linear equations in polynomial rings}\label{sec:flat sec}
Much of our work is done with Gr\"{o}bner bases, but it is sometimes convenient to use an older tool from Hermann's paper \cite{hermann1998question}.  

\begin{theorem}\label{thm:flat bd}
Let $k$ be a field and let $N_1,N_2\in \mathbb{N}$.  Consider the system of homogeneous equations $\{\sum_{j\leq N_2} f_{ij}X_j=0 \}_{i\leq N_1}$ where the coefficients $f_{ij}$ are polynomials in $k[x_1,\ldots,x_n]$ of total degree at most $d$. Then every solution in $(k[x_1,\ldots, x_n])^{N_2}$  of the system may be written as a linear combination (over $k[x_1,\ldots, x_n]$) of solutions that each have total degree at most $(2N_1d)^{2^n}$. 
\end{theorem}

\begin{proof}
See Theorem 3.2 of \cite{MR2051617} for this version of the bound; other proofs of the basic result are Theorem 1 of \cite{seidenberg1974constructions} and Theorem 2 of \cite{hermann1998question} for the original.
\end{proof}
\begin{remark}
In \cite{simmons_towsner}, Theorem 2.5, we call this fact ``internal flatness'' because it is equivalent to flatness of extensions of polynomial rings by internal polynomial rings in the sense of nonstandard analysis. \cite{MR2051617} refers to such results as ``effective flatness''.

\end{remark}

The preceding theorem generalizes to inhomogenous equations; here the relationship is with faithfully flat extensions.
\begin{theorem}\label{thm:fflat bd}
Let $k$ be a field and let $N_1, N_2\in \mathbb{N}$. Consider the system of  equations $\{\sum_{j\leq N_2} f_{ij}Y_j=h_i \}_{i\leq N_1}$, where the $f_{ij},h_i\in k[x_1,\ldots,x_n]$ have total degree at most $d$. If the equation has a solution in $(k[x_1,\ldots, x_n])^{N_2}$, then there exists a solution $y =(y_1,\ldots, y_{N_2})$ such that each $y_j$ has total degree at most $(2N_1d)^{2^n}$.
\end{theorem}

\begin{proof}
See Theorem 3.4 of \cite{MR2051617} or p.92 of \cite{renschuch1980beitrage} for this bound. Seidenberg \cite{seidenberg1974constructions} and Hermann  \cite{hermann1998question} gave earlier proofs but contain some errors (see p. 1,  \cite{MR2051617}).

\end{proof}

We also need a version of faithful flatness over the function field $k(x_1,\ldots,x_n)$.
\begin{theorem}\label{thm:L fflat}
Let $k$ be a field and let $N_1, N_2\in \mathbb{N}$. Consider the system of  equations $\{\sum_{j\leq N_2} f_{ij}Z_j=h_i \}_{i\leq N_1}$, where the $f_{ij},h_i\in k[x_1,\ldots,x_n]$ have total degree at most $d$. If the system has a solution in $(k(x_1,\ldots, x_n))^{N_2}$, then there exists a solution $z=(z_1,\ldots, z_{N_2}) \in (k(x_1,\ldots, x_n))^{N_2}$ such that each $z_j$ may be written as a ratio of polynomials of total degree at most $(2N_1d)^{2^n}$.
\end{theorem}

\begin{proof}
 Let $\widetilde{z}=(\widetilde{z}_1,\ldots,\widetilde{z}_{N_2})\in (k(x_1,\ldots, x_n))^{N_2}$ be a solution to the system. Multiply every equation in the system by the product of all denominators of  $\widetilde{z}_1,\ldots, \widetilde{z}_{N_2}$. It follows that there is a solution in $(k[x_1,\ldots, x_n])^{N_2+1}$ to the homogeneous system $\{\sum_{j\leq N_2} f_{ij}\widetilde{Z}_j=h_i\widetilde{Z} \}_{i\leq N_1}$ such that $\widetilde{Z}\neq 0$. Hence Theorem \ref{thm:flat bd} implies that there is a solution in $(k[x_1,\ldots, x_n])^{N_2+1}$ with each entry having total degree at most $(2N_1d)^{2^n}$ and  $\widetilde{Z}\neq 0$. Divide each entry by the nonzero value of $\widetilde{Z}$ to obtain the desired solution $z$ to the original system.

\end{proof}

The remaining flatness results require the field $k$ to have characteristic $p>0$.
\begin{theorem}[Based on Lemma 2.11 of \cite{schmidt89}] \label{thm:flat bd pth power}
Let $p,m, N_1, N_2\in \mathbb{N}$, with $p$ prime, and let $k$ be a field of characteristic $p$.  Consider the system of homogeneous equations $\{\sum_{j\leq N_2} f_{ij}Y_j^{p^m}=0 \}_{i\leq N_1}$ where the $f_{ij}$ are polynomials in $k[x_1,\ldots,x_n]$ of total degree at most $B$. Then every solution in $(k[x_1,\ldots, x_n])^{N_2}$  of the system may be written as a linear combination (over $k[x_1,\ldots, x_n]$) of solutions that each have total degree at most $(2N_1^2 N_2 \binom{B+n}{n}p^{m(n-1)}B)^{2^n}$. 
\end{theorem}

\begin{proof}
Let $V_0$ be a basis in $k$ of the set of $k^{p^m}$-linear combinations of the coefficients of all the $f_{ij}$. Define $V_1$ to be the set of all power products $\{x_1^{l_1}\cdot \ldots \cdot x_n^{l_n}\mid 0\leq l_j < p^m\}$. Finally, let $V=\{v_0v_1\mid v_0\in V_0, v_1\in V_1\}$. Since $V_0$ is a basis and each variable in an element of $V_1$ has power less than $p^m$, it follows that $V$ is linearly independent over $k^{p^m}[x_1^{p^m},\ldots,x_n^{p^m}]$. Each $f_{ij}$ may therefore be written uniquely as a sum $\sum_{v\in V} f_{ij,v}v $ where $f_{ij,v}\in k^{p^m}[x_1^{p^m},\ldots,x_n^{p^m}]$.

Substituting for $f_{ij}$ and changing the order of addition, the system $\{\sum_{j\leq N_2} f_{ij}Y_j^{p^m}=0\}_{i\leq N_1}$ becomes $\{\sum_{v\in V}\left(\sum_{j\leq N_2}f_{ij,v}Y_j^{p^m}\right)v=0\}_{i\leq N_1}.$ By linear independence, we may replace this with an equivalent system $\{\sum_{j\leq N_2} f_{ij,v} Y_j^{p^m}=0\}_{v\in V,i\leq N_1}$ whose coefficients belong to $k^{p^m}[x_1^{p^m},\ldots,x_n^{p^m}]$. This is a polynomial ring over the field $k^{p^m}$, so by \ref{thm:flat bd} every solution of $\{\sum_{j\leq N_2} f_{ij,v} Z_j=0\}_{v\in V, i\leq N_1}$ in $(k^{p^m}[x_1^{p^m},\ldots,x_n^{p^m}])^{N_2}$ is a $k^{p^m}[x_1^{p^m},\ldots,x_n^{p^m}]$-linear combination of solutions having bounded degree. That bound depends on the number of equations in the system and on the total degrees of  $f_{ij,v}$ as polynomials in $x_1^{p^m},\ldots,x_n^{p^m}$. The former is $N_1 |V|$, where $|V|$ is the cardinality of $V$. In turn, $|V|$ is the product of the cardinalities $|V_0|$ and $|V_1|$.  Lemma \ref{thm:num of power products} implies that $|V_0|\leq N_1N_2(\text{the number of monomials in any } f_{ij})\leq N_1N_2 \binom{B+n}{n}$. Note that $|V_1|\leq p^{mn}$, so $|V|=|V_0||V_1|\leq  N_1 N_2 \binom{B+n}{n}p^{mn}$. Since the $f_{ij,v}$ are obtained by factoring out powers of $x_1^{p^m},\ldots,x_n^{p^m}$ from monomials of the $f_{ij}$, the degree of the $f_{ij,v}$ as polynomials in $k^{p^m}[x_1^{p^m},\ldots,x_n^{p^m}]$ is at most $\frac{B}{p^m}$. Thus by \ref{thm:flat bd}, the desired bound is $(2 N_1|V|\frac{B}{p^m})^{2^n}\leq 
(2N_1^2 N_2 \binom{B+n}{n}p^{m(n-1)}B)^{2^n}$.

\end{proof}

\begin{theorem}\label{thm:baby fflat bd pth power}
Let $p,m, N_1, N_2\in \mathbb{N}$, with $p$ prime, and let $k$ be a field of characteristic $p$.   Consider the system of equations $\{\sum_{j\leq N_2} f_{ij}Y_j^{p^m}=h_i\}_{i\leq N_1} $ where $f_{ij},h_i$ are polynomials in $k[x_1,\ldots,x_n]$ of total degree at most $B$. If the system has a solution in $(k(x_1,\ldots, x_n))^{N_2}$, then there is a solution $y =(y_1,\ldots, y_{N_2}) \in (k(x_1,\ldots,x_n))^{N_2}$ such that each  $y_j$ may be written as a ratio of polynomials of total degree at most $(2N_1^2 (N_2+1) \binom{B+n}{n}p^{m(n-1)}B)^{2^n}$. 
\end{theorem}

\begin{proof}
We use the same approach as in Theorem \ref{thm:L fflat}. Consider the homogeneous system $\{\sum_{j\leq N_2} f_{ij}Y_j^{p^m}-h_iY^{p^m}=0\}_{i\leq N_1}$, where $Y$ is a new indeterminate. By \ref{thm:flat bd pth power}, every solution in $(k[x_1,\ldots, x_n])^{N_2+1}$  of the homogeneous system may be written as a $k[x_1,\ldots, x_n]$-linear combination of solutions that each have total degree at most $(2N_1^2 (N_2+1) \binom{B+n}{n}p^{m(n-1)}B)^{2^n}$. There is a solution in $(k(x_1,\ldots, x_n))^{N_2}$ to the original system, so by clearing denominators we obtain a solution in $(k[x_1,\ldots, x_n])^{N_2+1}$ to the homogeneous system for which $Y\neq 0$. Therefore there is a bounded solution $\widetilde{y}=(\widetilde{y}_1,\ldots, \widetilde{y}_{N_2},\widetilde{y}_{N_2+1})\in (k[x_1,\ldots, x_n])^{N_2+1}$ such that $\widetilde{y}_{N_2+1} \neq 0$. Divide $\widetilde{y}$ by $\widetilde{y}_{N_2+1}$ to obtain a bounded solution $(\widetilde{y}_1/\widetilde{y}_{N_2+1}, \ldots, \widetilde{y}_{N_2}/\widetilde{y}_{N_2+1},1)$ in $(k(x_1,\ldots,x_n))^{N_2+1}$. Then $y:=(\widetilde{y}_1/\widetilde{y}_{N_2+1},\ldots, \widetilde{y}_{N_2}/\widetilde{y}_{N_2+1})$ is the desired solution to the original equation.

\end{proof}
\subsection{Gr\"{o}bner Bases}

The basic properties of Gr\"{o}bner bases are laid out in many places; see, for instance, \cite{adams1994introduction,buchberger1998grobner,becker1993grobner, clo1,martin2000computational}. We need the following notions and facts. 

\begin{definition}
In the polynomial ring $k[x_1,\ldots,x_n]$, a \emph{power product} of variables $x_1,\ldots,x_n$ is a product $x_1^{r_1}\ldots x_n^{r_n}$ for some nonnegative integers $r_i$ (if all $r_i$ are zero, we write $1$ for the product). A \emph{monomial ordering} $<$ on 
$k[x_1,\ldots,x_n]$ is a well-ordering of the set of power products such that 
	\begin{enumerate}
	\item the product $1$ is the least element with respect to $<$ and  
	\item $<$ respects multiplication by power products: if $x_1^{r_1}\ldots x_n^{r_n} < x_1^{s_1}\ldots x_n^{s_n}$, then $(x_1^{t_1}\ldots x_n^{t_n})(x_1^{r_1}\ldots x_n^{r_n})=x_1^{r_1+t_1}\ldots x_n^{r_n+t_n}<x_1^{s_1+t_1}\ldots x_n^{s_n+t_n}=(x_1^{t_1}\ldots x_n^{t_n})(x_1^{s_1}\ldots x_n^{s_n})$. 

\end{enumerate}
\end{definition}

\begin{lemma}\label{thm:num of power products}
For positive integers $n,d$, there are $\binom{d+n-1}{n-1}$ power products of total degree $d$ in the indeterminates $x_1,\ldots,x_n$. Counting $1$ as a power product with each variable having degree 0, there are $\binom{d+n}{n}$ power products of total degree at most $d$.
\end{lemma}
\begin{proof}
This follows from an elementary counting argument; see, e.g., the ``stars and bars'' method in \cite{stanley1997enumerative}.
\end{proof}

\begin{definition}
Let $<$ be a monomial ordering on $k[x_1,\ldots,x_n]$. The \emph{leading term} $LT(f)$ of a nonzero polynomial $f$ is the monomial $ax^{r_1}\ldots x^{r_n}$ of $f$ such that the power product $x^{r_1}\ldots x^{r_n}$ is maximal with respect to $<$. The coefficient $a$ of the leading term is the \emph{leading coefficient} $LC(f)$ of $f$. (Note that a monomial is simply a power product multiplied by a field element.)
\end{definition}

\begin{definition}
Let $<$ be a monomial ordering on $k[x_1,\ldots,x_n]$ and let $I\subseteq k[x_1,\ldots,x_n]$ be an ideal. A set of nonzero polynomials $\{g_1,\ldots,g_t\}$ belonging to $I$ is a \emph{Gr\"{o}bner basis} (with respect to $<$) for $I$ if for every $f\in I$ there exists some $i$ such that the leading term $LT(g_i)$ of $g_i$ divides the leading term $LT(f)$ of $f$. 
\end{definition}

\begin{definition}
Let $F=\{f_1,\ldots,f_s\}\subseteq k[x_1,\ldots,x_n]$ be a finite set of nonzero polynomials. A polynomial $f\in k[x_1,\ldots,x_n]$ is \emph{reducible} with respect to $F$ if $f$ contains a nonzero monomial that is divisible by $LT(f_i)$ for some $1\leq i\leq s$. Otherwise we say that $f$ is \emph{reduced modulo} $F$. A \emph{reduced \GB basis} is one such that each generator is reduced with respect to the set containing all the others.
\end{definition}

\begin{theorem}
(Division algorithm for $k[x_1,\ldots,x_n]$)Let $F=\{f_1,\ldots,f_s\}\subseteq k[x_1,\ldots,x_n]$ be a finite set of nonzero polynomials and let $<$ be a monomial ordering on $k[x_1,\ldots,x_n]$. For every polynomial $f$, there exist $u_1,\ldots,u_s, r \in k[x_1,\ldots,x_n]$ such that 

	\begin{enumerate}
	\item $f=u_1f_1+\ldots+u_sf_s +r$,
	\item $r$ is reduced modulo $F$,
	\item  $LT(f)\geq LT(u_if_i)$  for all $1\leq i\leq s$, and
	\item  either $r$ is zero, or $f$ is reduced with respect to $F$ (in which case $u_i=0$ and $f=r$), or $LT(r)\leq LT(f)$.
        
	\end{enumerate} 
\end{theorem}

\begin{proof}
See Theorem 1.5.9 of \cite{adams1994introduction}.
\end{proof}

When $f=u_1f_1+\ldots+u_sf_s+r$ as in the theorem, we say $f$ \emph{reduces} to $r$ modulo $F$.

\begin{theorem}
Let $G$ be a \GB basis for a nonzero ideal $I$. A polynomial $f\in k[x_1,\ldots,x_n]$ belongs to $I$ if and only $f$ reduces to 0 modulo $G$. 
\end{theorem}

\begin{proof}
See Theorem 1.6.2 of \cite{adams1994introduction}.
\end{proof}

\begin{definition}
Let $f,g\in k[x_1,\ldots,x_n]$. Let $LCM(t_1, t_2)$ denote the least common multiple of two power products $t_1,t_2$. We define the $S$-polynomial $S(f,g)$ of $f,g$ to be $\displaystyle{\left(\frac{LCM(LT(f),LT(g))}{LT(f)}\right )\cdot f - \left(\frac{LCM(LT(f),LT(g))}{LT(g)}\right)\cdot g}$.   
\end{definition}

\begin{theorem}\label{thm:Buchberger}(Buchberger's Criterion)
Given a monomial ordering, let $G\subseteq k[x_1,\ldots,x_n]$ be a finite set of nonzero polynomials. $G$ is a \GB basis of the ideal generated by $G$ if and only if $S(f,g)$ reduces to 0 modulo $G$ for all $f,g\in G$.
\end{theorem}

\begin{proof}
See Theorem 1.7.4 \cite{adams1994introduction}.
\end{proof}

\begin{theorem}\label{thm:GB Dube}
Let $I$ be a nonzero ideal of $k[x_1,\ldots,x_n]$ generated by polynomials of total degree at most $d$. Then for any monomial ordering there is a reduced Gr\"{o}bner basis of $I$ whose elements have total degree at most $b_1(n,d):=\text{ min}\left\{2d^{2^{n}},2\left(\frac{d^2}{2} +d\right)^{2^{n-1}}\right\}$.  



\end{theorem}

\begin{proof}
The bound $2\left(\frac{d^2}{2} +d\right)^{2^{n-1}}$ was given by Dub\'{e} (see Corollary 8.3 of \cite{dube1990structure}). Note also that $2d^{2^n}$ is greater than or equal to Dub\'{e}'s bound whenever $d>1$. 
For $d=1$, the constant bound 1 suffices: it follows from Buchberger's criterion \ref{thm:Buchberger} that after Gaussian elimination we are left with a \GB basis (see Theorem 5.68 in \cite{becker1993grobner}). 


\end{proof}

\begin{remark}
We consider the more generous bound $2d^{2^n}$ because later we compose bounds with each other and we desire a relatively clean final answer.
\end{remark}


\begin{theorem}\label{thm:GB count} 
Let $I\subseteq k[x_1,\ldots,x_n]$ be an ideal generated by polynomials of total degree at most $d$. If $G=\{g_1,\ldots, g_t\}$ is a reduced \GB basis of $I$ whose elements have degree at most $b_1(n,d)$, then $t\leq (b_1(n,d)+1)^{n-1}$.
\end{theorem}
\begin{proof}

Since $G$ is reduced, no leading term of $g_i$ divides a leading term of $g_j$ for $i\neq j$. It suffices to find an upper bound on sets of monomials of total degree at most $b_1(n,d)$ such that no monomial in the set divides another. We induct on the number of variables. The constant value 1 clearly works for $n=1$.  Assume the claim holds for $n-1$.  For any $0\leq \alpha\leq b_1(n,d)$, consider the set $G_{x_1,\alpha}$ of leading terms of $G$ such that $x_1$ appears with degree $\alpha$. Since no element of $G_{x_1,\alpha}$ divides another, the same is true of the set $G_{x_1,\alpha}/(x_1^\alpha)$ of monomials formed by dividing the elements of $G_{x_1,\alpha}$ by $x_1^\alpha$. By the inductive hypothesis, $G_{x_1,\alpha}/(x_1^\alpha)$ (and hence $G_{x_1,\alpha}$) has at most $(b_1(n,d)+1)^{n-2}$ elements. $G$ is partitioned into $b_1(n,d)+1$ sets of size at most $(b_1(n,d)+1)^{n-2}$, so $G$ has at most $(b_1(n,d)+1)^{n-1}$ elements, as desired.
\end{proof}

\begin{remark}
As noted earlier, we aim for simplicity rather than sharpness. For a more detailed analysis of the number of generators in a \GB basis, see \cite{robbiano1990bounds}. The issue is closely related to Dickson's Lemma; see \cite{MR2858898,MR3448164} for bounds in that setting.
\end{remark}

We now convert the basic bound $b_1(n,d)$ into bounds on the auxiliary ideals that appear in the next section.




\begin{definition}
A monomial ordering $<$ of $k[y_1,\ldots y_m, x_1,\ldots,x_n]$ is an \emph{elimination ordering} if every power product in $k[x_1, \ldots x_n]$ (other than 1) is greater than every power product in $k[y_1,\ldots, y_m]$. If $I\subseteq k[y_1,\ldots y_m, x_1,\ldots,x_n]$ is an ideal, then $I\cap k[y_1,\ldots, y_m]$ is an \emph{elimination ideal} that eliminates the $x$-variables.
\end{definition} 

\begin{theorem}\label{thm:GB elim}
Let $<$ be an elimination ordering of  $k[y_1,\ldots y_m, x_1,\ldots,x_n]$ such that power products in the $x$-variables are greater than products in the $y$-variables. Let $G$ be a \GB basis with respect to $<$ of an ideal $I\subseteq k[y_1,\ldots y_m, x_1,\ldots,x_n]$ such that the elimination ideal $I\cap k[y_1,\ldots, y_m]$ is nonzero.  Then $G\cap k[y_1,\ldots,y_m]$ is a \GB basis of $I\cap k[y_1,\ldots, y_m]$. 
\end{theorem}
\begin{proof}
See Theorem 2.3.4 of \cite{adams1994introduction}.
\end{proof}

\begin{theorem}\label{thm:GB elim bd}

\smallskip
\noindent Let $<$ be an elimination ordering of  $k[y_1,\ldots y_m, x_1,\ldots,x_n]$ such that power products in the $x$-variables are greater than products in the $y$-variables. Let $I\subseteq k[y_1,\ldots y_m, x_1,\ldots,x_n]$ be generated by polynomials of total degree at most $d$. Also suppose that the elimination ideal $I\cap k[y_1,\ldots, y_m]$ is nonzero.  Then there there is a Gr\"{o}bner basis (with respect to $<$) of $I\cap k[y_1,\ldots,y_m]$ whose elements have total degree at most $b_1(m+n, d)$. 

\end{theorem}

\begin{proof}
By Theorem \ref{thm:GB elim}, the elimination ideal $I\cap k[y_1,\ldots,y_m]$ has a Gr\"{o}bner basis that is a subset of the \GB basis of $I$ given by \ref{thm:GB Dube}.

\end{proof}

\begin{theorem}\label{thm:bound for min polys}
Let $I\subseteq k[x_1,\ldots,x_n]$ be generated by polynomials of total degree at most $d$ and let $u\in k[x_1,\ldots,x_n]$ also have total degree at most $d$. If there exists a nonzero polynomial $p\in k[x_1,\ldots, x_r][Y]$ (where $Y$ is a new indeterminate and $r\leq n$) such that $p(u)\in I$, then there exists a nonzero polynomial $q\in  k[x_1,\ldots, x_r][Y]$ of total degree at most $b_1(n+1,d)$ such that $q(u)\in I$.
\end{theorem}
\begin{proof}
Consider the ideal $\widetilde{I}= (I, Y-u)\subseteq k[x_1,\ldots, x_n, Y]$. Since $p(u)\in I$, it follows that $p(Y)\in \widetilde{I}$. Use an elimination ordering in which power products of $x_{r+1},\ldots,x_n$ are greater than power products in $x_1,\ldots, x_r, Y$. Since the elimination ideal $\widetilde{I}\cap k[x_1,\ldots,x_r,Y]$ is nonzero, by Theorem \ref{thm:GB elim bd}, there is a Gr\"{o}bner basis of $\widetilde{I}\cap k[x_1,\ldots,x_r,Y]$ whose elements have total degree at most $b_1(n+1,d)$. Let $q$ be an element of this Gr\"{o}bner basis; note that $q(Y)\in \widetilde{I}$ implies that $q(u)\in I$.


\end{proof}

\subsection{Quotients and Saturations}
\begin{definition}
Let $I$ and $J$ be ideals of a commutative ring $R$. The \emph{quotient ideal} of $I$ by $J$  (denoted $I:J$) is defined to be $\{f\in R \mid Jf\subseteq I\}=\{f\in R\mid gf \in I \text{ for all } g \in J\}$.
\end{definition}

\begin{theorem} \label{thm:quotient}
\begin{enumerate} \item Let $I$ and $J$ be ideals of $k[x_1,\ldots, x_n]$ and let $w$ be a new variable. Then $I\cap J=(wI, (1-w)J)\cap k[x_1,\ldots,x_n]$.
\item Let $I$ be an ideal of $k[x_1,\ldots,x_n]$ let $f$ be a nonzero polynomial. Then the quotient ideal $I:(f)=\frac{1}{f}(I\cap(f)).$
\end{enumerate}

\end{theorem}
\begin{proof}
See Proposition 2.3.5 and Lemma 2.3.11 of \cite{adams1994introduction}.
\end{proof} 

\begin{theorem}\label{thm:GB quotient bd} 
If $I\subseteq k[x_1,\ldots, x_n]$  is a nonzero ideal generated by polynomials of total degree at most $d$ and $f$ is a nonzero polynomial of total degree at most $d$, then the quotient ideal $I:(f)$ has a  Gr\"{o}bner basis whose elements have total degree at most $b_1(n+1,d+1)$.

\end{theorem}
\begin{proof} 

By properties of a monomial ordering, the leading term of a product is the product of the leading terms. It follows from this and part 2 of \ref{thm:quotient} that by dividing the elements of a \GB basis for $I\cap(f)$ by $f$, we obtain a \GB basis for $I:(f)$. Hence it suffices to give a bound on \GB bases of $I\cap(f)$. By part 1 of \ref{thm:quotient}, $I\cap(f)= (wI, (1-w)f)\cap k[x_1,\ldots,x_n])$ where $w$ is a new variable. Use an elimination ordering with $w$ the greatest variable to eliminate $w$ from an ideal generated by polynomials of total degree at most $d+1$. Then \ref{thm:GB Dube} and \ref{thm:GB elim bd} imply that $I\cap(f)$ has a \GB basis of total degree at most $b_1(n+1,d+1)$. 

\end{proof}


\begin{definition}
Let $A$ be a commutative ring with 1. A subset $S$ of $A$ is a \emph{multiplicative set} if $1\in S, 0\notin S$, and $xy\in S$ for all $x,y\in S$. If $I$ is an ideal of $A$, the \emph{saturation} of $I$ by $S$ is the set of all elements $a\in A$ such that for some $s\in S$ we have $sa\in I$. We write $I:S^\infty$ to denote the saturation of $I$ by $S$.  If $b$ is a nonnilpotent element of $A$, we denote the multiplicative set $\{1, b, b^2, \ldots\}$ by $b^\infty$ and write the saturation as $I:b^\infty$.
\end{definition}

\begin{remark}
It is easy to see that $I:S^\infty$ is an ideal containing $I$. An equivalent characterization is that the saturation is the intersection $(S^{-1}I) \cap A$, where $S^{-1}I$ is the ideal generated by $I$ in the localization $S^{-1}A$. 
\end{remark}

In \ref{thm:base bd} we saturate an ideal $I\subseteq k[x_1,\ldots, x_n]$ by the multiplicative set $k[x_i]\setminus\! \{0\}$. 
To compute bounds on generators of this saturation we must  consider \GB bases  of polynomials with coefficients from more general rings. The definition of a \GB basis is slightly different from the field case to account for the fact that coefficients might be zero-divisors or fail to be units. (Monomial orderings on the set of power products are the same, however. Also, in this paper we only consider the case when the coefficient ring is a Noetherian integral domain.) For the precise definition, see Definition 4.1.13 of \cite{adams1994introduction}. The details are not necessary for our purposes because of the following theorem:  
\begin{theorem}\label{thm:GB field to dom}
Let $G$ be a \GB basis (in the sense of fields) of an ideal $I\subseteq k[y, x_1,\ldots, x_n]$ with respect to an elimination order with the $x$-variables greater than $y$. Then $G$ is a \GB basis (in the sense of Noetherian integral domains) of $I$ viewed as an ideal of $(k[y])[x_1,\ldots,x_n]$. 
\end{theorem} 

\begin{proof}
See Theorem 4.1.18 in \cite{adams1994introduction}.
\end{proof}

\begin{theorem}\label{thm:sat by g} 
Let $R$ be a Noetherian integral domain, let $g\in R[y_1,\ldots,y_m]$ be a nonzero polynomial, and let $I$ be an ideal of $R[y_1,\ldots,y_m]$. Letting $w$ be a new variable, we have $I:g^\infty=(I, wg-1)\cap R[y_1,\ldots, y_m]$.
\end{theorem}

\begin{proof} 
See Proposition 4.4.1 in \cite{adams1994introduction}.
\end{proof}

\begin{theorem}\label{thm:sat by R} 
Let $R$ be a Noetherian integral domain and let $I\subseteq R[y_1,\ldots,y_m]$ be an ideal having a \GB basis $G=\{g_1,\ldots, g_t\}$ (in the sense of Noetherian integral domains). Let $g=LC(g_1)LC(g_2)\cdots LC(g_t)$ be the product of the leading coefficients of the elements of $G$. Then $I:(R\setminus\!\{0\})^\infty=I:g^\infty$. 
\end{theorem}
\begin{proof}
See Proposition 4.4.4 in \cite{adams1994introduction}. 
\end{proof}

\begin{theorem}\label{thm:GB sat bd} 
If $I\subseteq k[x_1,\ldots, x_n]$  is a nonzero ideal generated by polynomials of total degree at most $d$, then the saturation $I:(k[x_1]\setminus\!\{0\})^\infty$ has a  Gr\"{o}bner basis whose elements have total degree at most $b_2(n,d):=b_1(n+1,(b_1(n,d)+1)^{n-1}b_1(n,d) +1)$. 

\end{theorem}
\begin{proof}

Choose an elimination ordering having $x_1$ as least variable. By \ref{thm:GB Dube} and \ref{thm:GB count}, $I$ has a reduced \GB basis $G=\{g_1,\ldots,g_t\}$ with respect to this ordering such that the elements of $G$ have total degree at most $b_1(n,d)$ and $t\leq (b_1(n,d)+1)^{n-1}$.
Now view $I$ as an ideal of $(k[x_1])[x_2,\ldots,x_n])$, the polynomial ring in $n-1$ variables over $k[x_1]$. By Theorem \ref{thm:GB field to dom} we know that $G$ is still a \GB basis in the sense of coefficient rings that are Noetherian integral domains. By Theorem \ref{thm:sat by R}, $I:(k[x_1]\setminus\!\{0\})^{\infty}= I:g^\infty$, where $g=LC(g_1)LC(g_2)\cdots LC(g_t)$ is the product of the leading coefficients of the elements of $G$. Thus the degree of $g$ in $x_1$ is bounded by $tb_1(n,d)$. 

Theorem \ref{thm:sat by g} shows that $I:g^\infty=(I, wg-1)\cap (k[x_1])[x_2,\ldots, x_n]$, where $w$ is a new variable. Hence $I:(k[x_1]\setminus\!\{0\})^\infty=(I, wg-1)\cap k[x_1, x_2,\ldots, x_n]$ and we may view $(I, wg-1)$ as an ideal of $k[x_1,x_2,\ldots, x_n,w]$ generated by polynomials of total degree at most $tb_1(n,d)+1$.

Use an elimination ordering with $w$ the greatest variable to eliminate $w$ from $(I, wg-1)$. By \ref{thm:GB Dube}, \ref{thm:GB elim bd}, and the bound on $t$, the resulting elimination ideal has a \GB basis of total degree at most $b_1(n+1,tb_1(n,d)+1)\leq b_1(n+1,(b_1(n,d)+1)^{n-1}b_1(n,d) +1)$.  

\end{proof}

\begin{theorem}
$b_2(n,d)\leq 2^{n2^{n+3}}d^{n2^{2n+1}}$
\end{theorem}
\begin{proof}
We get a bound of the desired form through the following chain of inequalities:

\begin{align*} &b_2(n,d)=b_1(n+1,(b_1(n,d)+1)^{n-1}b_1(n,d) +1)\leq 2((b_1(n,d)+1)^{n-1}b_1(n,d) +1)^{2^{n+1}} \\ & \leq 2((2d^{2^n}+1)^{n-1}\cdot 2d^{2^{n}}+1)^{2^{n+1}} \leq 
2((2^2d^{2^n})^{n-1}\cdot 2^2d^{2^n})^{2^{n+1}} =2(2^{2n}d^{n2^n})^{2^{n+1}} \\ & \leq 2^{n2^{n+3}}d^{n2^{2n+1}}.
\end{align*}




\end{proof}

\begin{section}{Proof of the main theorem}\label{sec:main theorem}
We continue to use $k$ to denote an arbitrary field. $I$ denotes a proper ideal generated by polynomials $f_1,\ldots, f_s\in k[x_1,\ldots, x_n]$ of total degree at most $d$. As explained in the outline, we choose a maximal set of indeterminates that are algebraically independent modulo $I$.  Relabeling if necessary, we may assume they are $x_1,\ldots,x_r$. Since 0 is the only polynomial in $k[x_1,\ldots x_r]\cap I$, whenever $fg\in I$ for $f\in k[x_1,\ldots,x_r]\setminus \{0\}$, either $g\in I$ or $f,g$ are counterexamples to primality of $I$.   

\begin{definition}\label{def:bdd prime and max}

If $I\subseteq k[x_1,\ldots,x_n]$ is a proper ideal and $b$ is a natural number, we say that $I$ is \emph{prime up to $b$} if for all $f,g\in k[x_1,\ldots,x_n]$ of total degree at most $b$ such that $fg\in I$, either $f\in I$ or $g\in I$. By \emph{maximal up to $b$}, we mean that for every $f\not\in I$ with degree $\leq b$, $f$ has an inverse mod $I$: there is some $g$ with $fg-1\in I$. 
\end{definition}

Before proceeding, we take care of the simple cases in the main theorem:

\begin{theorem}\label{thm:easy cases} 
Let $I\subseteq k[x_1,\ldots,x_n]$ be a proper ideal generated by polynomials of degree at most $d$. Let $\{x_1,\ldots,x_r\}$ be a maximal subset of $\{x_1,\ldots,x_r,\ldots, x_n\}$ that is algebraically independent modulo $I$.
\begin{itemize}
\item  If $n=1$, then $I$ is prime if $I$ is prime up to $d$.
\item  If $d=1$ (with $n$ arbitrary), then $I$ is prime.
\item If $r=0$, then $I$ is prime if $I$ is prime up to $nb_1(n,d)$.
\item If $r=n$, then $I=\{0\}$ and is prime.
\end{itemize}
\end{theorem}

\begin{proof}

First let $n=1$. $k[x_1]$ is a principal ideal domain and $I$ is of the form $(f)$ where $f$ is the greatest common divisor of the given generators of $I$. $f$ has degree at most $d$ and $I$ is prime if and only if $f$ is irreducible (if $f$ is reducible, no proper factor belongs to $I$). If $I$ is prime up to $d$, then $f$ is irreducible.

When $d=1$, so $I$ is generated by linear polynomials, $I$ is necessarily prime; this follows geometrically from the irreducibility of affine linear varieties. Therefore primality up to any natural number (e.g., 0) vacuously guarantees primality if $d=1$.

If $r=0$, by Theorem \ref{thm:GB elim bd} there is a nonzero polynomial $\alpha_i(x_i)\in I\cap k[x_i]$ of degree at most $b_1(n,d)$ for each $x_1,\ldots, x_n$. Using $\alpha_i$ to reduce higher powers of $x_i$ mod $I$, we see that any $f$ is equivalent mod $I$ to some $\widetilde{f}$ of total degree at most $nb_1(n,d)$. Hence $fg\in I$ if and only if $\widetilde{f}\widetilde{g}\in I$ and then primality up to $nb_1(n,d)$ ensures primality.

The final claim is clear.
\end{proof}

\begin{remark}
 If any $x_i\in I$, then modulo $I$ we may assume that $x_i$ does not appear in any polynomials we consider.  Our claimed bound on detecting prime ideals increases with $n$, so having fewer variables only makes the bound easier to prove.
\end{remark}

Due to the preceding theorem and remark, we may safely assume in the following work that $n,d>1$, $0<r<n$ and that $x_i\notin I$ for any $1\leq i \leq n$.

Denote by $J$ the ideal generated by $I$ in $k(x_1,\ldots,x_r)[x_{r+1},\ldots,x_n]$. The key properties of $J$ are that if $I$ is prime up to a sufficient level, then $I=J\cap k[x_1,\ldots, x_n]$ and $J$ is maximal. Together these prove that $I$ is prime. 


\begin{subsection}{Proof of $I=J\cap k[x_1,\ldots, x_n]$} 
To obtain this result we rely heavily on the bounds from Section \ref{prelim} on Gr\"{o}bner bases for various kinds of ideals.


\begin{theorem}\label{thm:base bd}
Let $1\leq i \leq n$. If $f\in k[x_i]$ is irreducible and there exists $g\in k[x_1,\ldots,x_n]$ such that $fg\in I$ but $g\notin I$, then the degree of $f$ is at most $b_3(n,d):=b_1(n+1,b_2(n,d)+1)$.

\end{theorem}
\begin{proof}

Let $q_1,\ldots, q_N$ be generators of the saturation ideal $I:(k[x_i]\setminus\{0\})^\infty\subseteq k[x_1,\ldots,x_n]$. By \ref{thm:GB sat bd}, the total degree of each generator is at most $b_2(n,d)$. For each $1\leq l\leq N$, consider the ideal $(I:(q_l))\cap k[x_i]$. Since $k[x_i]$ is a PID, each ideal $(I:(q_l))\cap k[x_i]$ is generated by some polynomial $u_l\in k[x_i]$; in particular, $u_lq_l\in I$. By \ref{thm:GB quotient bd} and \ref{thm:GB elim bd}, the degree of $u_l$ is at most $b_3(n,d):=b_1(n+1,b_2(n,d)+1)$. Let $f\in k[x_i]$ be irreducible and let $g$ be an element of $k[x_1,\ldots,x_n]\setminus I$ with the property that $fg\in I$. We claim that $\text{deg}(f)\leq \text{ max}_{1\leq l\leq N}\{\text{deg}(u_l)\}\leq b_3(n,d)$. 

Note that $g\in I:(k[x_i]\setminus\{0\})^\infty$ and so there are $\alpha_l\in k[x_1,\ldots,x_n]$ such that $g=\sum_{l=1}^N \alpha_lq_l$. Multiplying both sides by $\prod_{l=1}^N u_l$, we obtain $\left(\prod_{l=1}^N u_l\right)g =\left(\prod_{l=1}^N u_l\right)\left(\sum_{l=1}^N \alpha_lq_l\right)=\sum_{l=1}^{N}\left(\prod_{l=1}^N u_l\right)\alpha_lq_l\in I$, so $\prod_{l=1}^N u_l\in  (I:(g))\cap k[x_i]$. Also note that $f\in (I:(g))\cap k[x_i]$ and that $(I:(g))\cap k[x_i]$ is proper since $g\notin I$. Given that $f$ belongs to the same  proper ideal of $k[x_i]$ as $\prod_{l=1}^N u_l$ and is irreducible, $f$ must generate the ideal and hence divide $u_l$ for some $l$. In particular, the degree of $f$ is at most the maximum of the degrees of the $u_l$ and is less than or equal to $b_3(n,d)$ as claimed.


\end{proof}

\begin{remark}
Note that the result holds for any variable $x_i$, not necessarily one of the $x_1,\ldots,x_r$. Also, the proof goes through unchanged for any field $L$ and any ideal of $L[x_1,\ldots,x_n]$ generated by elements of degree at most $d$. This is important because in \ref{thm:second ind step} we apply \ref{thm:base bd} after changing the field from $k$ to $k(x_{i_1},\ldots,x_{i_{l}})$ for some subset $\{x_{i_1},\ldots,x_{i_{l}}\}$ of the variables.
\end{remark}
\begin{theorem}
$b_3(n,d)\leq 2^{n2^{2n+6}}d^{n2^{3n+2}}$.
\end{theorem}

\begin{proof}
We compute as follows:

\begin{align*}
&b_3(n,d) =b_1(n+1,b_2(n,d)+1)\leq 2(2^{n2^{n+3}}d^{n2^{2n+1}}+1)^{2^{n+1}}   \\ & \leq 2(2^{n2^{n+4}}d^{n2^{2n+1}})^{2^{n+1}} =
2(2^{n2^{2n+5}}d^{n2^{3n+2}})\leq 2^{n2^{2n+6}}d^{n2^{3n+2}} .
\end{align*}
\end{proof}
We aim to show that $J\cap k[x_1,\ldots, x_n]=I$. We consider arbitrary subsets $\{x_{i_1},\ldots, x_{i_l}\}\subseteq\{x_1,\ldots,x_r\}$ of cardinality $l$ having certain properties and show inductively that $\{x_{i_1},\ldots, x_{i_l}, x_{i_{l+1}}\}$ retains the properties. Let $J_l$ denote $Ik(x_{i_1},\ldots, x_{i_l})[\text{remaining variables from }x_1,\ldots, x_n ]$ and let $J_{l+1}$ denote $Ik(x_{i_1},\ldots, x_{i_l},x_{i_{l+1}})[\text{remaining variables from }x_1,\ldots, x_n ]$. With this notation we have $J=J_r$.

\begin{theorem}\label{thm:first ind step} 
Let $1\leq l\leq r$ and suppose that $I:(f)=I$ for all irreducible $f\in k[x_{i_1},\ldots, x_{i_l}]$ and every subset $\{x_{i_1},\ldots, x_{i_l}\}\subseteq\{x_1,\ldots,x_r\}$ of cardinality $l$.  Then

\begin{enumerate}
\item $I:(f)=I$ for all $f\in k[x_{i_1},\ldots, x_{i_l}]\setminus\{0\}$ (irreducible or not) and
\item $J_l\cap k[x_1,\ldots,x_n]=I$.
\end{enumerate}
\end{theorem}

\begin{proof}
The reverse containment is trivial for both parts. For the forward direction, let $f$ factor as $ \prod_{l=1}^Np_l$ for some irreducible polynomials $p_l\in k[x_{i_1},\ldots, x_{i_l}]$ and suppose $fg\in I$. It follows that $\left(\prod_{l=2}^N p_l\right)g \in I:(p_1)$, which is equal to $I$ by the hypothesis since $p_1$ is irreducible. Continuing this way, we conclude that $g\in I$. This proves (1).

For the second claim, let $g\in J_l\cap k[x_1,\ldots,x_n]$. Clearing denominators, we see that $fg\in I$ for some $f \in  k[x_{i_1},\ldots, x_{i_l}]$. Part (1) implies that $g\in I$.
\end{proof}

\begin{theorem}\label{thm:second ind step} 
Let $1\leq l<r$ and assume $n>1$. Suppose that $J_l\cap k[x_1,\ldots,x_n]=I$ for all $\{x_{i_1},\ldots, x_{i_l}\}\subseteq\{x_1,\ldots,x_r\}$ and $I:(f)=I$  for all irreducible $f\in k[x_{i_1},\ldots, x_{i_l},x_{i_{l+1}}]$ of total degree at most $(n-1)\cdot b_3(n,d)$. Then $J_{l+1}\cap k[x_1,\ldots,x_n]=I$.
\end{theorem} 
\begin{proof}
By \ref{thm:first ind step}, it suffices to show that $I:(f)=I$ for all irreducible $f\in k[x_{i_1},\ldots, x_{i_l},x_{i_{l+1}}]$. Pick any such $f$. By hypothesis we already have $I:(f)=I$ if the total degree of $f$ is less than or equal to $(n-1)\cdot b_3(n,d)$, so suppose that the degree of $f$ exceeds this. At least one of the $l+1\leq r\leq n-1$ variables  must appear with degree greater than $b_3(n,d)$. Without loss of generality we may assume that variable is $x_{i_{l+1}}$ (because we assume the properties in the statement for arbitrary subsets of $\{x_1,\ldots,x_r\}$ of size $l$). By Gauss' Lemma, $f$ remains irreducible over $k(x_{i_1},\ldots, x_{i_l})$. Since \ref{thm:base bd} applies to irreducible univariate polynomials over any field and the $x_{i_{l+1}}$-degree of $f$ is greater than $b_3(n,d)$, we deduce that $J_{l}:(f)=J_{l}$. Let $g\in k[x_1,\ldots,x_n]$ be such that $fg\in I$. 
We conclude from $J_{l}:(f)=J_{l}$ that $g\in J_l$. It follows that $g\in I$ because by hypothesis $J_l\cap k[x_1,\ldots,x_n]=I$. This proves that $I:(f)=I$ as required.  

\end{proof}
\begin{theorem}\label{thm:bdd prime quotient}
Let $n,d>1$. If $I$ is prime up to $b_4(n,d):=b_1(n+1, (n-1)\cdot b_3(n,d)+1)$, then for all $1\leq l \leq r$, subsets $\{x_{i_1},\ldots, x_{i_l}\}\subseteq\{x_1,\ldots,x_r\}$, and irreducible $f\in k[x_{i_1},\ldots, x_{i_l}]$ of total degree at most $(n-1)\cdot b_3(n,d)$ we have $I:(f)=I$.
\end{theorem}

\begin{proof}
It suffices for the generators of $I:(f)$ to belong to $I$. By \ref{thm:GB quotient bd}, such a generator $g$ has total degree at most $b_4(n,d):=b_1(n+1, (n-1)\cdot b_3(n,d)+1)$.  Since $fg\in I$ but $f\notin I$, by primality up to $b_4(n,d)$ we conclude that $g\in I$.
\end{proof}

\begin{theorem}\label{thm:J intersect}
Let $n,d>1$. If $I$ is prime up to $b_4(n,d)$, then $J\cap k[x_1,\ldots,x_n] =I$.
\end{theorem} 
\begin{proof}
We induct on the size $l$ of subsets $\{x_{i_1},\ldots, x_{i_l}\}$ of $\{x_1,\ldots,x_r\}$ to prove that $J_l\cap k[x_1,\ldots,x_n]=I$ for each $1\leq l\leq r$. 
By \ref{thm:bdd prime quotient}, $I:(f)=I$ for all irreducible $f\in k[x_{i_1},\ldots, x_{i_l}]$ of total degree at most $(n-1)\cdot b_3(n,d)$. (This satisfies one of the two hypotheses of \ref{thm:second ind step}.) Letting $l=1$, Theorem \ref{thm:base bd} now implies that $I:(f)=I$ for irreducible $f\in k[x_l]$ of any degree. Theorem \ref{thm:first ind step} then yields the base case $J_1\cap k[x_1,\ldots, x_n]$. We can now apply \ref{thm:second ind step} repeatedly to obtain $J\cap k[x_1,\ldots,x_n] =I$. 

\end{proof}

\end{subsection}

\begin{subsection}{$J$ is maximal}

We start with an important lemma that allows us to control the degrees of $x_{r+1},\ldots,x_n$ modulo $I$:
\begin{lemma}\label{thm:deg of min polys}
Let $I$ be defined as before and suppose $I$ is prime up to $b_1(n,d)$. 
\begin{enumerate}
\item 
For each $r+1\leq j\leq n$, there is a polynomial $w_j \in I\cap  k[x_1,\ldots,x_r][x_j]$ such that $w_j$ has total degree at most $b_1(n,d)$, is irreducible, and remains irreducible over $k(x_1,\ldots,x_r)[x_j]$.
\item $k(x_1,\ldots, x_r)[x_{r+1},\ldots, x_n]/J$ is a finite-dimensional $k(x_1,\ldots,x_r)$-vector space of dimension at most  $\prod_{j=r+1}^n d_j\leq (b_1(n,d))^{n-r}$, where $d_j$ is the $x_j$-degree of $w_j$ in the first part of the lemma.
\end{enumerate}


\end{lemma}
\begin{proof}
By definition of $\{x_1,\ldots,x_r\}$ as a maximal algebraically independent set mod $I$, we know that $I\cap k[x_1,\ldots,x_r]=\{0\}$ and that for all $j>r$ there exist elements of positive $x_j$-degree in $I\cap k[x_1,\ldots,x_r][x_j]$. By Theorem \ref{thm:GB elim bd}, $I\cap k[x_1,\ldots, x_r][x_j]$ has a Gr\"{o}bner basis whose elements have total degree at most at most $b_1(n,d)$; the $x_j$ degree of each generator is positive. Consider the set $W=\{w\in I\cap k[x_1,\ldots,x_r][x_j]\mid \text{the total degree of } w \text{ is at most } b_1(n,d)\}$ and let $N$ be the minimal $x_j$-degree of any element of $W$. 
Let $W_N\subseteq W$ be the subset of $W$ whose elements have $x_j$-degree $N$. Choose $w_{j}=a_Nx_j^N +\ldots + a_0\in W_N$ such that $a_N\in k[x_1,\ldots,x_r]$ has minimal total degree out of all coefficients of $x_j^N$ in $W_N$.

We claim that $w_{j}$ is irreducible in $k[x_1,\ldots, x_r][x_j]$. To see this, suppose toward contradiction that $w_{j}$ properly factors as $fg=(b_lx_j^l +\ldots+ b_0)(c_mx_j^m+\ldots + c_0)$. Then either $b_l\in k[x_1,\ldots,x_r]$ would have lower total degree than $a_N$ or $l<N$, and likewise for $c_m$ and $m$. Since $f,g$ have total degree at most $b_1(n,d)$, 
by primality up to $b_1(n,d)$ we have $f\in I$ or $g\in I$.  But this contradicts either minimality of $N$ or that of $a_N$, so $w_{j}$ does not factor in $k[x_1,\ldots, x_r][x_j]$. By Gauss' lemma \ref{thm:Gauss lemma}, $w_{j}$ remains irreducible over $k(x_1,\ldots, x_r)[x_j]$.

The final claim follows from the existence of the $w_j$ and the definition of $d_j$.
\end{proof}

\begin{lemma}\label{thm:M exists}
$J$ is contained in a maximal ideal $M$ of $k(x_1,\ldots, x_r)[x_{r+1},\ldots, x_n]$.
\end{lemma}
\begin{proof}
We need only show that $J$ is proper. If on the contrary $1$ is a $k(x_1,\ldots,x_r)[x_{r+1},\ldots, x_n]$-linear combination of $f_1,\ldots, f_s$, then by clearing  denominators in an equation witnessing that $1\in J$ we obtain a nonzero element $k[x_1,\ldots,x_r]\cap I$. This contradicts algebraic independence of $x_1,\ldots, x_r$ mod $I$. 
\end{proof}

We prove that $J$ is maximal by showing that $J=M$. A crucial step is the following:

\begin{subsubsection}{Proof that elements of $M$ having small degree belong to $J$}
When working in the ring $k(x_1,\ldots,x_r)[x_{r+1},\ldots,x_n]$, we must keep track of the variables $x_1,\ldots, x_r$ as well as $x_{r+1},\ldots,x_n$. For convenience we use the following terminology:
\begin{definition}
Let $\alpha/\beta\in k(x_1,\ldots,x_r)\setminus \{0\}$ with $\alpha,\beta$ coprime as polynomials in $k[x_1,\ldots,x_r]$. We say the $(x_1,\ldots,x_r)$-\emph{degree} (or just degree) of $\alpha/\beta$ is the maximum of the total degrees of $\alpha$ and $\beta$. A monomial $\left(\alpha/\beta\right)T$, where $T$ is a power product in $x_{r+1},\ldots,x_n$, is defined to have $(x_1,\ldots,x_r)$-degree equal to the sum of the $(x_1,\ldots,x_r)$-degree of $\alpha/\beta$ and the total degree of $T$. Given $f\in k(x_1,\ldots,x_r)[x_{r+1},\ldots,x_n]\setminus\{0\}$, we define the $(x_1,\ldots,x_r)$-degree of $f$ to be the maximum of the $(x_1,\ldots,x_r)$-degrees of the monomials appearing in $f$.

Let $\mathcal{J}\subseteq k(x_1,\ldots,x_r)[x_{r+1},\ldots,x_n]$ be a proper ideal and $B$ a natural number. We say $\mathcal{J}$ is $(x_1,\ldots,x_r)$-\emph{prime} up to $B$ if for all $f,g\in k(x_1,\ldots,x_r)[x_{r+1},\ldots,x_n]$ of $(x_1,\ldots,x_r)$-degree at most $B$ such that $fg\in \mathcal{J}$, either $f\in \mathcal{J}$ or $g\in \mathcal{J}$. Likewise, if every $f\in k(x_1,\ldots,x_r)[x_{r+1},\ldots,x_n]\setminus \mathcal{J}$ of $(x_1,\ldots,x_r)$-degree at most $B$ is invertible mod $\mathcal{J}$, we say that $\mathcal{J}$ is $(x_1,\ldots,x_r)$-\emph{maximal} up to $B$. 
\end{definition}

\begin{lemma}\label{thm:clear denom}
Let $f\in k(x_1,\ldots,x_r)[x_{r+1},\ldots,x_n]\setminus\{0\}$ have $(x_1,\ldots,x_r)$-degree at most $B$. Then there exists a polynomial $g\in k[x_1,\ldots,x_r]$ of total degree at most $B\cdot \binom{B+n}{n}$ such that $gf\in k[x_1,\ldots,x_r, x_{r+1},\ldots,x_n]$ has total degree at most $B\left( 1+\binom{B+n}{n}\right) $.
\end{lemma}

\begin{proof}
Write $f$ as $\sum_i \left(\alpha_i/\beta_i\right)T_i$, where $\alpha_i,\beta_i$ are coprime polynomials in $k[x_1,\ldots,x_r]$ and $T_i$ is a power product in $x_{r+1},\ldots, x_n$. Since $f$ has $(x_1,\ldots,x_r)$-degree bounded by $B$, each $\alpha_i,\beta_i$, and $T_i$ has total degree at most $B$. Clear denominators and keep track of the degrees. By \ref{thm:num of power products}, there are $\binom{B+n}{n}$ power products of degree at most $B$ in $n$ variables, so $f$ contains at most $\binom{B+n}{n}$ monomials. Let $g$ be $\prod_i \beta_i$, which has total degree at most $B\cdot \binom{B+n}{n}$. It follows that $gf=\left(\prod_i \beta_i\right)\cdot\sum_i \left(\alpha_i/\beta_i\right)T_i\in k[x_1,\ldots,x_r, x_{r+1},\ldots,x_n]$ and has total degree at most $B\left( 1+\binom{B+n}{n}\right) $.
\end{proof}

\begin{lemma}\label{thm:prime to rational prime}
For any $d_0$, if $I$ is prime up to $\text{max}\left\{b_4(n,d), d_0\! \left( 1+\binom{d_0+n}{n}\right)\right\} $, then $J=Ik(x_1,\ldots,x_r)[x_{r+1},\ldots,x_n]$ is $(x_1,\ldots,x_r)$-prime up to $d_0$.
\end{lemma}
\begin{proof}

Let $g,h\in k(x_1,\ldots,x_r)[x_{r+1},\ldots,x_n]$ each have $(x_1,\ldots,x_r)$-degree $\leq d_0$ and suppose $gh\in J$. By \ref{thm:clear denom}, there are are $g_0,h_0\in  k[x_1,\ldots,x_r]$ such that $g_0g, h_0h \in  k[x_1,\ldots,x_n]$ and the total degrees of $g_0g, h_0h$ are at most $d_0\!\left( 1+\binom{d_0+n}{n}\right)$. By \ref{thm:J intersect},  $(g_0g)(h_0h)\in I$ since $(g_0g)(h_0h)\in J\cap k[x_1,\ldots,x_n] $ and $I$ is prime up to at least $b_4(n,d)$.  Since $I$ is also prime up to at least $ d_0\! \left( 1+\binom{d_0+n}{n}\right) $, either $g_0g\in I$ or $h_0h\in I$, showing that either $g$ or $h$ belongs to $J$.


\end{proof}

\begin{lemma} [Based on {\cite{schmidt89}}, Lemma 2.3]\label{thm:schmidt_2.3}


Let $f\in k(x_1,\ldots,x_r)[x_{r+1},\dots,x_n]$ and let the $(x_1,\ldots,x_r)$-degree of $f$ be at most $d_0$, where $d\leq d_0$. 
  For each $j=r+1,\dots, n$, denote by $d_j$ the $x_j$-degree of the polynomial $w_j$ obtained in \ref{thm:deg of min polys}. There exists $\widetilde{f}\in k(x_1,\ldots,x_r)[x_{r+1},\dots,x_n]$ such that the following hold:

\begin{enumerate}
\item $\text{deg}_{x_j}\widetilde{f} < d_j$,
\item $f-\widetilde{f} \in J=Ik(x_1,\ldots,x_r)[x_{r+1},\dots, x_n]$, and
\item the $(x_1,\ldots,x_r)$-degree of $\widetilde{f}$ is bounded by $d_0(1+b_1(n,d))^{n-r}$.

\end{enumerate}

\end{lemma}

\begin{proof}

The claim follows from successively dividing $f$ by the polynomials $w_{r+1},\dots, w_{n}$. By Lemma \ref{thm:deg of min polys}, each such $w_j \in I\cap  k[x_1,\ldots,x_r][x_j]\subseteq J$ has total degree at most $b_1(n,d)$, so $d_j\leq b_1(n,d).$ The first division step involving $w_{r+1}$ at most adds the $(x_1,\ldots,x_r)$-degrees of $f$ and $w_{r+1}$ while the $x_{r+1}$-degree decreases. (If $\text{deg}_{x_{r+1}}f$ is already less than $d_{r+1}$, we may skip to $w_{r+2}$.) At each step we increase the original $(x_1,\ldots,x_r)$-degree of $f$ (which was at most $d_0$) by at most $b_1(n,d)$, and we can continue at most $d_0$-many times. Therefore, the $(x_1,\ldots,x_r)$-degree of the remainder after completing  division by $w_{r+1}$ is bounded by $d_0 + d_0\cdot b_1(n,d)= d_0(1+b_1(n,d))$. 

After dividing by $w_j$, the degree of the remainder in $x_j$ is less than $d_j$ and can never increase thereafter because $x_j$ does not appear in $w_l$ for $r<l\neq j$. The remainder $\widetilde{f}$ after considering all $n-r$ polynomials $w_j$ thus satisfies the first two claims. The pattern repeats, giving us a bound of $d_0 + d_0\cdot \binom{n-r}{1}\cdot b_1(n,d)+\ldots+ d_0\cdot \binom{n-r}{n-r-1}\cdot(b_1(n,d))^{n-r-1}+ d_0\cdot (b_1(n,d))^{n-r}=d_0(1+b_1(n,d))^{n-r}$ when the process terminates. 

\end{proof}

\begin{lemma}[Based on Lemma 2.4 of \cite{schmidt89}]\label{thm:lem2.4}
  Let $f\in k(x_1,\ldots,x_r)[x_{r+1},\ldots,x_n]$ with $(x_1,\ldots,x_r)$-degree $\leq d_0$ where $d\leq d_0$.  For each $j=r+1,\dots, n$, denote by $d_j$ the $x_j$-degree of the polynomial $w_j$ obtained in \ref{thm:deg of min polys}. There is a monic polynomial $\theta_f\in k(x_1,\ldots,x_r)[Y]$ such that:
  \begin{enumerate}
  \item the $Y$-degree of $\theta_f$ is at most $\prod_{j=r+1}^n d_j$,
\item $\theta_f(f)\in J=Ik(x_1,\ldots,x_r)[x_{r+1},\ldots,x_n]$, and
  \item the $(x_1,\ldots,x_r)$-degree of $\theta_f$ is at most $(\prod_{j=r+1}^n d_j)(d_0(1+b_1(n,d))^{n-r})$.
  \end{enumerate}
\end{lemma}
\begin{proof}
   Let $\{v_1,\ldots,v_t\}$ be power products in $x_{r+1},\ldots, x_n$ that form a basis of $k(x_1,\ldots,x_r)[x_{r+1},\ldots,x_n]$ mod $J$. By Lemma \ref{thm:deg of min polys}, $t$ is at most $\prod_{j=r+1}^n d_j$.   By Lemma \ref{thm:schmidt_2.3}, 
  for each $v_i$ there is a $\widetilde{fv_i}$ such that:
  \begin{itemize}
  \item for $r<j\leq n$, the $x_j$-degree of $\widetilde{fv_i}$ is less than $d_j$,
  \item $fv_i-\widetilde{fv_i}\in J$, and
\item $\widetilde{fv_i}$ has $(x_1,\ldots,x_r)$-degree at most $d_0(1+b_1(n,d))^{n-r}$.
  \end{itemize}
So we may write
\[\widetilde{fv_i}=\sum_{l=1}^ta_{il}v_l\]
where each $a_{il}\in k(x_1,\ldots,x_r)$ has $(x_1,\ldots,x_r)$-degree at most $d_0(1+b_1(n,d))^{n-r}$.

We use the determinant trick. Consider the $t\times t$-matrix $A(Y)=[\delta_{il}Y-a_{il}]$ and 
note that the determinant of $A(Y)$ is monic and $A(f)$ must belong to $J$:  if $v$ is the column vector with entries $v_l$, then $[a_{il}]v =fv$ mod $J$. Thus $[\delta_{il} f -a_{il}]v =0$ mod $J$. Multiply both sides by the adjugate of $[\delta_{il} f -a_{il}]$ to obtain $[\delta_{il}\,\mathrm{det}(A(f))]v \in J$. The power product 1 appears among the $v_i$, so $\mathrm{det} (A(f))\in J$.  

Let $\theta_f=\mathrm{det}(A(Y))$.  Then $\theta_f$ has $Y$-degree at most $\prod_{j=r+1}^n d_j$ and $(x_1,\ldots,x_r)$-degree $\leq t\cdot(d_0(1+b_1(n,d))^{n-r})\leq (\prod_{j=r+1}^n d_j)(d_0(1+b_1(n,d))^{n-r}).$  
\end{proof}

\begin{theorem}[Based on Lemma 2.7 of \cite{schmidt89}] \label{thm:weak bounded maximal}
Let $d_0\in \mathbb{N}$ and let $\widetilde{d}= (\prod_{j=r+1}^n d_j)(d_0)(1+(1+b_1(n,d))^{n-r})$, with $d_j$ defined in Lemma \ref{thm:deg of min polys}.  If $I$ is prime up to $\text{max}\left\{b_4(n,d), \widetilde{d}\! \left( 1+\binom{\widetilde{d}+n}{n}\right)\right\} $, then $J=Ik(x_1,\ldots,x_r)[x_{r+1},\ldots,x_n]$ is $(x_1,\ldots,x_r)$-maximal up to $d_0$.



\end{theorem}
\begin{proof}
Suppose we have $f\in k(x_1,\ldots,x_r)[x_{r+1},\ldots,x_n]\setminus J$ with $(x_1,\ldots,x_r)$-degree $\leq d_0$.  By Lemma \ref{thm:lem2.4}, there is a monic $\theta_f\in k(x_1,\ldots,x_n)[Y]$ with $Y$-degree $\leq \prod_{j=r+1}^n d_j$ and $(x_1,\ldots,x_r)$-degree at most $(\prod_{j=r+1}^n d_j)(d_0(1+b_1(n,d))^{n-r})$ so that $\theta_f(f)\in J$. Note that $\theta_f(f)$ has $(x_1,\ldots,x_r)$-degree at most $(d_0)(\prod_{j=r+1}^n d_j) + (\prod_{j=r+1}^n d_j)(d_0(1+b_1(n,d))^{n-r})=(\prod_{j=r+1}^n d_j)(d_0)(1 +(1+b_1(n,d))^{n-r})=\widetilde{d}$. 

Let us write $\theta_f=Y^m\theta'$ with $m$ maximal, so $f^m\theta'(f)\in J$.  By Lemma \ref{thm:prime to rational prime}, primality of $I$ up to $\text{max}\left\{b_4(n,d), \widetilde{d}\! \left( 1+\binom{\widetilde{d}+n}{n}\right)\right\} $ implies  $(x_1,\ldots,x_r)$-primality of $J$ up to $\widetilde{d}$ .  Since $f^m\theta'(f)$ has $(x_1,\ldots,x_r)$-degree $\leq \widetilde{d}$ and $f\not\in J$, we must have $\theta'(f)\in J$. This also implies that the constant term of $\theta'(Y)$ cannot be 0.

Since $\theta'(Y)=Y^s+\sum_{i=1}^{s-1}a_iY^i+a_0$ with $a_0\neq 0$ and each $a_i\in k(x_1,\ldots,x_r)$, we have
\[(a_0^{-1}f^{s-1}+\sum_{i=1}^{s-1}a_0^{-1}a_if^{i-1})f+1\in J,\]
and therefore $f$ is invertible mod $J$.
\end{proof}

\begin{cor}\label{thm:M down to J}
Let $d_0, \widetilde{d}, I, J,M$ be defined as in Theorem \ref{thm:weak bounded maximal}, with $I$ being prime up to $\text{max}\left\{b_4(n,d), \widetilde{d}\! \left( 1+\binom{\widetilde{d}+n}{n}\right)\right\} $. If $f\in M\subseteq k(x_1,\ldots,x_r)[x_{r+1},\ldots,x_n]$ has total $(x_1,\ldots,x_r)$-degree less than or equal to $d_0$, then $f\in J$.
\end{cor}
\begin{proof}
If $f\notin J$, then $f$ is invertible mod $J$ because $J$ is $(x_1,\ldots, x_r)$-maximal up to $d_0$ by Theorem \ref{thm:weak bounded maximal}. But since $f\in M$, it follows that $1\in M$, contradicting that $M$ is a proper ideal.
\end{proof}

\end{subsubsection}

\begin{subsubsection}{Proof that $J=M$}
 We argue that $J=M$ using induction on the number of variables $x_{r+1},\ldots, x_n$. It is convenient to first treat the case that $x_{r+2}+M,\ldots, x_n+M$ are separable over $k(x_1,\ldots, x_r)$ (technically, the isomorphic copy $k(x_1,\ldots, x_r)+M$ contained in the field $k(x_1,\ldots, x_r)[x_{r+1},\ldots,x_n]/M$). Theorem \ref{thm:separable max induction}, which is based on Theorem 2.8 in \cite{schmidt89}, 
contains the heart of the matter.

\begin{lemma}\label{thm:subfields mod M}
Let $M$ be defined as above and let $M_j:=M\cap k(x_1,\ldots,x_r)[x_{r+1},\ldots, x_j]$ for $r+1\leq j\leq n$. Then $M_j$ is a maximal ideal of $k(x_1,\ldots,x_r)[x_{r+1},\ldots, x_j]$ and the field $k(x_1,\ldots, x_r)[x_{r+1},\ldots,x_j]/M_j$ embeds in $k(x_1,\ldots, x_r)[x_{r+1},\ldots,x_n]/M$.
\end{lemma}
\begin{proof}
Consider the generators $x_{r+1}+M,\ldots,x_n+M$ of the field extension $k(x_1,\ldots, x_r)[x_{r+1},\ldots,x_n]/M$ over $k(x_1,\ldots, x_r)$. Because each generator is algebraic, the $k(x_1,\ldots, x_r)$-algebra $k(x_1,\ldots, x_r)[x_{r+1}+M,\ldots,x_j+M]$ is a subfield of $k(x_1,\ldots, x_r)[x_{r+1},\ldots,x_n]/M$. The $k(x_1,\ldots, x_r)$-algebra is isomorphic to $k(x_1,\ldots, x_r)[x_{r+1},\ldots,x_j]/M_j$ via the obvious map sending $x_{r+1} +M_j$ to $x_{r+1}+M$, etc. 
\end{proof}

\begin{theorem}[Based on Theorem 2.8 in \cite{schmidt89}]\label{thm:separable max induction}
Let $I,J,M$ be defined as above. Suppose that $I$ is prime up to $b_5(n,d):=(2d)^{2^{3n^2+2n}}$ and  $x_{r+2}+M,\ldots, x_n+M$ are separable over $k(x_1,\ldots, x_r)$. Then for each $r+1\leq j \leq n$ there exist $c_{r+1,j}, \ldots, c_{j,j} \in k[x_1,\ldots, x_r], h_j \in k[x_1,\ldots, x_r][Y]$, and  $\phi_{r+1,j},\ldots, \phi_{j,j} \in k(x_1,\ldots, x_r)[Y]$ such that, if we define $U_j:=c_{r+1,j}x_{r+1} + c_{r+2,j}x_{r+2}+\ldots + c_{j,j}x_j$ and $M_j:=M\cap k(x_1,\ldots,x_r)[x_{r+1},\ldots, x_j]$,

\begin{enumerate}
\item $U_j+ M_j$ generates $k(x_1,\ldots, x_r)[x_{r+1},\ldots,x_j]/M_j$ over $k(x_1,\ldots,x_r)$,
\item $h_j$ is the minimal polynomial of $U_j+M_j$ over $k(x_1,\ldots, x_r)$,
\item $h_j(U_j)\in I$ and $x_{r+1} - \phi_{r+1,j}(U_j),\ldots, x_j- \phi_{j,j}(U_j)$ all belong to $J$, 
\item and the total degrees of $c_{r+1,j}, \ldots, c_{j,j} \in k[x_1,\ldots, x_r]$ and of $h_j \in k[x_1,\ldots, x_r][Y]$ are at most $b_5(n,d)$.
\end{enumerate}

\end{theorem}



\begin{proof}

We induct on $r+1\leq j \leq n$. Let $w_j(x_j)$ be as defined in \ref{thm:deg of min polys}. For the base case, define $c_{r+1,r+1}=1$, $h_{r+1} = w_{r+1}(Y)$, and $\phi_{r+1,r+1}=x_{r+1}-w_{r+1}(Y)$. Clearly $U_{r+1}+ M_{r+1}= x_{r+1}+ M_{r+1}$ generates $k(x_1,\ldots, x_r)[x_{r+1}]/M_{r+1}$ over $k(x_1,\ldots,x_r)$. Lemma \ref{thm:deg of min polys} implies that claims 2 and 4 hold for $c_{r+1,r+1}$ and $ h_{r+1}$. The bounds suffice since the total degree of $w_{r+1}$ is at most $b_1(n,d)\leq 2d^{2^n}\leq b_5(n,d)$.   Claim 3 holds because $x_{r+1} - \phi_{r+1,r+1}(U_{r+1})=w_{r+1}(x_{r+1})\in I\subseteq J$. 

Suppose 1-4 hold for $U_j,c_{r+1,j},\ldots, c_{j,j},  h_{j}$, and $\phi_{r+1,j},\ldots, \phi_{j,j}$. We may also assume that the total degrees of $U_j$ and $h_j$ are at most $(2d)^{2^{3jn+2j}}$ (this holds for the base case $j=r+1$).  Set $c_{r+1,j+1}, \ldots, c_{j,j+1}$ equal to $c_{r+1,j}, \ldots, c_{j,j}$, respectively. By \ref{thm:subfields mod M}, $k(x_1,\ldots, x_r)[x_{r+1},\ldots,x_{j+1}]/M_{j+1}$ is a field extending (an isomorphic copy of) $k(x_1,\ldots, x_r)[x_{r+1},\ldots,x_j]/M_j$. From the inductive hypothesis, $h_j$ is the minimal polynomial of $U_j+M_{j+1}$ over $k(x_1,\ldots, x_r)$. 
By assumption $x_{j+1}+M_{j+1}$ is separable and by Lemma \ref{thm:deg of min polys} $w_{j+1}$ is the minimal polynomial of $x_{j+1} + M_{j+1}$ over $k(x_1,\ldots, x_r)$.  We also know that $k(x_1,\ldots, x_r)$ is an infinite field and $U_{j}+M_j $ generates $k(x_1,\ldots, x_r)[x_{r+1}\ldots,x_j]/M_j$. Hence the primitive element theorem \ref{thm:Primitive element} implies that there exists $c_{j+1,j+1}\in k[x_1,\ldots, x_r]$ such that $c_{j+1,j+1}$ has total degree at most $(\text{deg }h_j)(b_1(n,d)-1)\leq (2d)^{2^{3jn+2j}}\cdot (2d^{2^n}-1)$ and $U_{j+1}+ M_{j+1}=U_j + c_{j+1,j+1}x_{j+1}+ M_{j+1}$ generates $k(x_1,\ldots, x_r)[x_{r+1},\ldots,x_{j+1}]/M_{j+1}$ over $k(x_1,\ldots,x_r)$. 
 Note that the total degree of $U_{j+1}$ is at most $(2d)^{2^{3jn+2j}}\cdot (2d^{2^n}-1)+1\leq (2d)^{2^{3jn+2j}}\cdot 2d^{2^n}\leq (2d)^{2^{3jn+2j} +2^n}$. 

Because $k(x_1,\ldots, x_r)[x_{r+1},\ldots, x_n]/J$ is a finite-dimensional $k(x_1,\ldots,x_r)$-vector space, $U_{j+1}$ is algebraic mod $J$ (and hence mod $I$). 
By Theorem \ref{thm:bound for min polys}, the set $H=\{h\in k[x_1,\ldots,x_r][Y]\setminus\{0\}\mid \text{  the total degree of } h \text{ is at most } b_1(n+1, \text{max}\{d,\text{deg } U_{j+1}\}) \text{ and }h(U_{j+1})\in I\}$ is nonempty. The total degrees of elements of $H$ are bounded by  

\begin{align*}
&b_1(n+1, \text{max}\{d,\text{deg } U_{j+1}\})\leq b_1(n+1, (2d)^{2^{3jn+2j} +2^n})\leq 2((2d)^{2^{3jn+2j} +2^n})^{2^{n+1}} \\&
 \leq (2d)^{(2^{3jn+2j} +2^n)2^{n+1} +1} \leq (2d)^{2^{3jn+2j+n+1} +2^{2n+1} +1} \leq (2d)^{2^{3(j+1)n+2(j+1)}}. 
\end{align*}

\medskip
By the same argument we used in Lemma \ref{thm:deg of min polys}, there exists some irreducible $h_{j+1}\in H $; by definition of $H$, the total degree of $h_{j+1}$ is at most $(2d)^{2^{3(j+1)n+2(j+1)}}\leq b_5(n,d)$.  (Since $j\leq n$, primality up to $b_5(n,d):=(2d)^{2^{3n^2+2n}}$ suffices for each $j$. The bound $(2d)^{2^{3(j+1)n+2(j+1)}}$ has the correct form to continue the induction.)

It follows that $h_{j+1}$ is the minimal polynomial of $U_{j+1}+M_{j+1}$ over $k(x_1,\ldots, x_r)$. Note that $k(x_1,\ldots, x_r)[x_{r+1},\ldots,x_{j+1}]/M_{j+1}$ is isomorphic to $k(x_1,\ldots, x_r)[Y]/(h_{j+1})$ under the map sending $U_{j+1}+M_{j+1}$ to $Y+(h_{j+1})$.

It remains to define $\phi_{r+1,j+1},\ldots, \phi_{j+1,j+1}$. We start with $\phi_{j+1,j+1}$ and use it to determine the others. Consider the polynomial ring $k(x_1,\ldots, x_r)[Y,Z]$. The proof of the primitive element theorem  \ref{thm:Primitive element} shows 
that there exist $A,B\in k(x_1,\ldots, x_r)[Y,Z]$ such that $Z-(x_{j+1}+M_{j+1})= h_j((U_{j+1}+ M_{j+1})-c_{j+1,j+1}Z)A(U_{j+1}+M_{j+1}, Z) + w_{j+1}(Z)B(U_{j+1}+M_{j+1}, Z)$. (This uses the fact that $h_j$ and $w_{j+1}$ are the minimal polynomials of $U_j +M_{j+1}= U_{j+1}-c_{j+1,j+1}x_{j+1} + M_{j+1}$ and $x_{j+1}+M_{j+1}$, respectively.)

Since $U_{j+1}$ generates $k(x_1,\ldots, x_r)[x_{r+1},\ldots,x_{j+1}]$ mod $M_{j+1}$, there is some polynomial $C\in k(x_1,\ldots, x_r)[Y]$ such that $x_{j+1}-C(U_{j+1}) \in M_{j+1}$. It follows that $Z-(C(U_{j+1})+M_{j+1})= h_j((U_{j+1}+ M_{j+1})-c_{j+1,j+1}Z)A(U_{j+1}+M_{j+1}, Z) + w_{j+1}(Z)B(U_{j+1}+M_{j+1}, Z)$. By the isomorphism sending $U_{j+1}+M_{j+1}$ to $Y+(h_{j+1})$ we have $Z-C(Y)= h_j(Y-c_{j+1,j+1}Z)A(Y, Z) + w_{j+1}(Z)B(Y, Z) + h_{j+1}(Y)D(Y,Z)$ for some $D\in k(x_1,\ldots, x_r)[Y,Z]$. Define $\phi_{j+1,j+1}(Y)$ to be $C(Y)$ and observe that upon substituting $x_{j+1}$ for $Z$ and $U_{j+1}$ for $Y$ we get $x_{j+1} -\phi_{j+1,j+1}(U_{j+1})\in J$.

By hypothesis $x_{r+1}-\phi_{r+1,j}(U_j),\ldots, x_j-\phi_{j,j}(U_j) \in J$. Since $U_j = U_{j+1}-c_{j+1,j+1}x_{j+1}$ and $x_{j+1} -\phi_{j+1,j+1}(U_{j+1})\in J$, we define $\phi_{l,j+1}(Y)$ to be $\phi_{l,j}(Y -c_{j+1,j+1}\phi_{j+1,j+1}(Y))$ for $r+1\leq l \leq j$. This implies that $x_l-\phi_{l,j+1}(U_{j+1})\in J$, completing the proof.

\end{proof}
\begin{theorem}\label{thm:separable maximal}
Let $I,J,M$ be defined as above. Suppose $I$ is prime up to $b_5(n,d)$ and  $x_{r+2}+M,\ldots, x_n+M$ are separable over $k(x_1,\ldots, x_r)$. Then $J$ is a maximal ideal of $k(x_1,\ldots, x_r)[x_{r+1},\ldots,x_n]$.
\end{theorem}
\begin{proof}
We prove that $J=M$ by showing for any $p\in k(x_1,\ldots,x_r)[x_{r+1},\ldots,x_n]$ that if $p\notin J$, then $p\notin M$ ($J\subseteq M$ by definition). We use the polynomials $U_n, x_{r+1}- \phi_{r+1,n}(U_n)\ldots, x_n-\phi_{n,n}(U_n)$ and $h_n(U_n)$ from \ref{thm:separable max induction}.  By the theorem, $U_n+M_n=U_n+M$ generates the field $k(x_1,\ldots,x_r)[x_{r+1},\ldots,x_n]/M$, $x_{r+1}- \phi_{r+1,n}(U_n)\ldots, x_n-\phi_{n,n}(U_n)$ all belong to $J$, and $h_n(U_n)\in I$ with $h_n$ being the minimal polynomial of $U_n+M$. Hence we may replace, modulo $J$, each $x_{r+1},\ldots, x_n$ in $p$ with a polynomial in $U_n$; i.e., $p=p_1(U_n) + q$ for some $q\in J$ and some $p_1\in k(x_1,\ldots,x_r)[Y]$. 
Since $p\notin J$ but $h_n(U_n)\in I$, we know that $p_1(U_n)\notin J$ and the remainder  $p_2(U_n)$ from dividing $p_1(U_n)$ by $h_n(U_n)$ is not zero. 

Because $p_2(U_n)$ is a nonzero polynomial of lower degree in $U_n$ than $h_n(U_n)$, minimality of $h_n$ implies that $p_2(U_n)\notin M$. Since $p_1(U_n)$ is equal to $p_2(U_n)$ mod $J$, it follows that $p\notin M$.
\end{proof}

We now reduce the inseparable case to the separable. Theorem \ref{thm:maximal} is based on Theorems 2.8 and 2.12 in \cite{schmidt89}. We use induction, with the base step depending on \ref{thm:separable maximal}. The inductive step requires a specialization of the faithful flatness results from Section \ref{sec:flat sec}. In Lemma \ref{thm:L fflat mod J} we continue to use the same field $k$ and ideals $I\subseteq k[x_1,\ldots, x_r, x_{r+1},\ldots, x_n]$ and $J\subseteq k(x_1,\ldots, x_r)[x_{r+1},\ldots, x_n]$ as before. 

\begin{lemma}[Based on Lemma 2.10 of \cite{schmidt89}]\label{thm:L fflat mod J}
Let $g\in k[x_1,\ldots, x_{r}, \ldots, x_n]$ have total degree at most $B\geq d$ and let $D=\{t_1,\ldots, t_N\}$ be a set of power products $x_{r+1}^{i_{r+1}}\cdot \ldots \cdot (x_n)^{i_n}$, with each element of $D$ having total degree at most $B$. Suppose there exist $\beta_t \in k(x_1,\ldots,x_r)$ for each $t\in D$ such that $g + \sum_{t\in D}\beta_t t \in J$. Then there exist $\gamma_t\in k(x_1,\ldots, x_r)$ such that $g + \sum_{t\in D}\gamma_t t \in J$ and each $\gamma_t$ may be written as a ratio of polynomials of total degree at most $\left( 2\binom{(2B)^{2^n}+d + n-r}{n-r} B \right)^{2^r}$.
\end{lemma}

\begin{proof}
Considering $g$ and the generators $f_1,\ldots, f_s$ of $I$ as polynomials in $k(x_1,\ldots, x_r)[x_{r+1},\ldots, x_n]$, we may apply Theorem \ref{thm:fflat bd} to conclude that there are $g_i\in k(x_1,\ldots, x_r)[x_{r+1},\ldots, x_n]$ having total degree (in $x_{r+1},\ldots,x_n$) at most $(2B)^{2^n}$ such that  $g + \sum_{t\in D}\beta_t t = \sum_{i=1}^s g_if_i$. The monomials in $x_{r+1},\ldots, x_n$ that appear in the equation have total degree at most $(2B)^{2^n}+d$. By Lemma \ref{thm:num of power products}, there are $N_1:=\binom{(2B)^{2^n}+d + n-r}{n-r}$ power products in $x_{r+1},\ldots, x_n$ of total degree at most $(2B)^{2^n}+d$. 

From the $k(x_1,\ldots, x_r)$-coefficients of the power products of $x_{r+1},\ldots, x_n$ in $g + \sum_{t\in D}\beta_t t = \sum_{i=1}^s g_if_i$, we obtain a system of $N_1$-many linear equations over $k[x_1,\ldots, x_r]$ whose coefficients have degree at most $B$. Theorem \ref{thm:L fflat} implies that the system has solutions in $k(x_1,\ldots,x_r)$ bounded by $\left( 2\binom{(2B)^{2^n}+d + n-r}{n-r}B \right)^{2^r}$. Thus there are $\gamma_t,g_i'\in k(x_1,\ldots,x_r)$ that satisfy the claimed bound such that $g + \sum_{t\in D}\gamma_t t = \sum_{i=1}^s g_i'f_i\in J$.


\end{proof}

We are ready to handle the general case. To keep track of the bounds, we use the following names in Theorem \ref{thm:maximal}. While $b_1(n,d), b_4(n,d)$, and $b_5(n,d)$ are the same as before, note that the other values are not necessarily identical to those in the statements of earlier theorems.

\begin{notation}\label{def:intermed bds}
\begin{align*}
N_1&:=\left( 2\binom{(2(n-r)b_1(n,d))^{2^n}+d + n-r}{n-r}(n-r)b_1(n,d) \right)^{2^r},\\
N_2&:=\left( 2\binom{(2(n-r)(b_1(n,d))^2)^{2^n}+d + n-r}{n-r}(n-r)(b_1(n,d))^2 \right)^{2^r},\\
N_3&:=b_1(n,d)^{n-r},\\
N_4&:=N_3N_2 +N_1,\\
N_5&:=(2N_3^2 (N_3+1) \binom{N_4+r}{r}(b_1(n,d))^{(r-1)}N_4)^{2^r},\\
N_6&:=N_5+(n-r)(b_1(n,d))^2,\\
\widetilde{d}&:=N_3N_6(1+(1+b_1(n,d))^{n-r}),\\
b_6(n,d)&:= max\left\{b_4(n,d), b_5(n,d), \widetilde{d}\! \left( 1+\binom{\widetilde{d}+n}{n}\right)\right\}.
\end{align*}
\end{notation}

\begin{theorem}[Based on Theorems 2.8 and 2.12 of \cite{schmidt89}]\label{thm:maximal}
If $I$ is prime up to $b_6(n,d)$, then $J$ is a maximal ideal of $k(x_1,\ldots, x_r)[x_{r+1},\ldots,x_n]$.
\end{theorem}

\begin{proof}
Since $I$ is prime at least up to $b_5(n,d)$, by Theorem \ref{thm:separable maximal} we are done unless $k$ has prime characteristic $p$ and at least some of $x_{r+1}, \ldots, x_n$ are inseparable mod $M$. For each $r+1\leq j \leq n$, let $h_j(Y)\in k(x_1,\ldots,x_r)[Y]$ be the minimal polynomial of $x_j$ mod $M$. 
If $h_j$ happens to be separable, define $m_j$ to be 0. If $h_j$ is inseparable, then every exponent of $Y$ that appears in $h_j$ is divisible by $p$. 
Let $m_j$ be the greatest power of $p$ such that $p^{m_j}$ divides all exponents of $Y$ that appear in $h_j$. 
Recall that the polynomial $w_j(x_j)\in k[x_1,\ldots,x_r][x_j]$ from Lemma \ref{thm:deg of min polys} belongs to $I\subseteq M$. Being a minimal polynomial, $h_j(Y)$ divides $w_j(Y)$ in $k(x_1,\ldots,x_r)[Y]$. Lemma \ref{thm:deg of min polys} then implies that $p^{m_j}\leq \text{ deg}_Y(h_j)\leq \text{ deg}_Y(w_j)=d_j \leq b_1(n,d).$

Let $\widetilde{h}_j(Z) \in k(x_1,\ldots,x_r)[Z]$ be the polynomial obtained from $h_j$ by replacing $Y^{p^{m_j}}$ with $Z$. Note that $\widetilde{h}_j$ remains irreducible. It follows from our choice of $m_j$ that not every exponent of $Z$ appearing in $\widetilde{h}_j$ is divisible by $p$. Hence $\widetilde{h}_j$ is the minimal polynomial of $x_j^{p^{m_j}}$ mod $M$ and $x_j^{p^{m_j}}+M$ is separable over $k(x_1,\ldots,x_r)$.

Let $m$ be the greatest of $m_{r+1},\ldots, m_n$ and consider the polynomial ring $k(x_1,\ldots,x_r)[x_{r+1}^{p^{m}},\ldots, x_n^{p^{m}}]$. Observe that $M\cap k(x_1,\ldots,x_r)[x_{r+1}^{p^{m}},\ldots, x_n^{p^{m}}]$ is still maximal: if $f\in k(x_1,\ldots,x_r)[x_{r+1}^{p^{m}},\ldots, x_n^{p^{m}}]\setminus M$, then by maximality of $M$, there exists $g\in 
k(x_1,\ldots,x_r)[x_{r+1},\ldots,x_n]$ such that $fg-1\in M$. Since $k$ has prime characteristic $p$, $(fg-1)^{p^m}=f^{p^m}g^{p^m}-1= f(f^{p^m -1}g^{p^m}) -1\in M$, with  $f^{p^m -1}g^{p^m}\in k(x_1,\ldots,x_r)[x_{r+1}^{p^{m}},\ldots, x_n^{p^{m}}]$. This shows that $M\cap k(x_1,\ldots,x_r)[x_{r+1}^{p^{m}},\ldots, x_n^{p^{m}}]$ is maximal.
Since $x_{r+1}^{p^{m}}+M\cap k(x_1,\ldots,x_r)[x_{r+1}^{p^{m}},\ldots, x_n^{p^{m}}],\ldots,x_n^{p^{m}}+M\cap k(x_1,\ldots,x_r)[x_{r+1}^{p^{m}},\ldots, x_n^{p^{m}}]$ are separable over $k(x_1,\ldots, x_r)$, by \ref{thm:separable maximal} the ideal $J\cap k(x_1,\ldots,x_r)[x_{r+1}^{p^{m}},\ldots, x_n^{p^{m}}]=M\cap k(x_1,\ldots,x_r)[x_{r+1}^{p^{m}},\ldots, x_n^{p^{m}}]$. We use induction to show that $J=M$. 

For $r\leq l\leq n$, define $S_{l}=J\cap k(x_1,\ldots,x_r)[x_{r+1}, \ldots, x_{l}, x_{l+1}^{p^{m}},\ldots,  x_n^{p^{m}}]$ and $M_l=M\cap k(x_1,\ldots,x_r)[x_{r+1}, \ldots, x_{l}, x_{l+1}^{p^{m}},\ldots, x_n^{p^{m}}]$. (So $M_{n}=M\cap k(x_1,\ldots,x_r)[x_{r+1}, \ldots, x_n]=M$.) As was the case with $l=r$,  each $M_l$ is a maximal ideal of $k(x_1,\ldots,x_r)[x_{r+1}, \ldots, x_{l}, x_{l+1}^{p^{m}},\ldots, x_n^{p^{m}}]$. 
Define the following quotient rings:

\begin{align*}
A_l=k(x_1,\ldots,x_r)[x_{r+1}, \ldots ,x_{l}, x_{l+1}^{p^{m}},\ldots, x_n^{p^{m}}]/S_{l}
\end{align*}

\noindent and 

\begin{align*}
B_l=k(x_1,\ldots,x_r)[x_{r+1}, \ldots ,x_{l},x_{l+1}^{p^{m}},\ldots, x_n^{p^{m}}]/M_l.
\end{align*}
\noindent The previous paragraph established that $S_{r}=M_{r}$, the base case of our induction. Now suppose that $S_j=M_j$ for some $r\leq j< n$. In particular,  $A_j=B_j$. We claim that $S_{j+1}=M_{j+1}$.



Note that $A_{j+1}$ surjects as a ring onto the field $B_{j+1}$ via the natural map sending $\alpha+S_{j+1}$ to $\alpha + M_{j+1}$ for $\alpha \in k(x_1,\ldots,x_r)[x_{r+1}, \ldots ,x_{j+1},x_{j+2}^{p^{m}},\ldots, x_n^{p^{m}}]$. We prove that this map is injective by showing that the vector space dimension of $B_{j+1}$ over $B_j$ is greater than or equal to that of $A_{j+1}$ over the field $A_j=B_j$. This will show that $A_{j+1}=B_{j+1}$ and $S_{j+1}=M_{j+1}$. 

Since $x_{j+1}^{p^m}+M_{j+1}\in B_j$ and $B_{j+1}$ is generated as a field over $B_j$ by $x_{j+1}+M_{j+1}$, $B_{j+1}/B_j$ is a purely inseparable extension. Hence the minimal polynomial of $x_{j+1}+M_{j+1}$ over $B_j$ has the form $X^{p^{\mu}} + b$ for some $\mu\in \mathbb{N}$ and $b\in B_j$. In particular, the degree of the extension $B_{j+1}/B_j$ is $p^{\mu}$. It remains to show that the vector space dimension of $A_{j+1}$ over $A_j=B_j$ is at most $p^{\mu}$.

Denote by $\tau$ the minimal natural number such that $x_{j+1}^{p^{\tau}}+ S_{j+1} \in A_j$. Note that $\mu\leq \tau\leq m$. Let $D=\{t_1+S_j, \ldots, t_N+S_j\}$ be a $k(x_1,\ldots,x_r)$-basis of $A_j=B_j$, where $t_1,\ldots, t_N$ are representative power products in the variables $x_{r+1}, \ldots x_j, x_{j+1}^{p^m},\ldots, x_n^{p^m}$. Thus there are  $\alpha_t, \beta_t\in k(x_1,\ldots,x_r)$ such that  $x_{j+1}^{p^{\mu}}+ \sum_{t\in D} \alpha_t t\in M$ and $x_{j+1}^{p^{\tau}}+ \sum_{t\in D} \beta_t t\in J$. (From now on we write $t\in D$ to refer to any of $t_1,\ldots, t_N$.) It follows from \ref{thm:deg of min polys} that we may assume that a power product $t\in D$ has total degree at most $(n-r)b_1(n,d)$ as a polynomial in $x_{r+1}, \ldots x_j, x_{j+1}^{p^m},\ldots, x_n^{p^m}$. We have $p^\tau\leq p^m\leq b_1(n,d)$, so by Lemma \ref{thm:L fflat mod J} we may assume that the $\beta_t$ are ratios of polynomials of total degree at most $N_1:=\left( 2\binom{(2(n-r)b_1(n,d))^{2^n}+d + n-r}{n-r}(n-r)b_1(n,d) \right)^{2^r}$.

If $\tau=\mu$, then $x_{j+1}^{p^{\tau}}+ \sum_{t\in D} \beta_t t\in J$ implies that the dimension of $A_{j+1}$ over $A_j$ is at most $p^\mu$. Hence it suffices to obtain a contradiction from the assumption $\mu<\tau$. Suppose then that $\mu<\tau$ and consider 

\begin{align}\label{block0}
(x_{j+1}^{p^{\mu}}+ \sum_{t\in D} \alpha_t t)^{p^{\tau-\mu}}-(x_{j+1}^{p^{\tau}}+ \sum_{t\in D} \beta_t t) &=\sum_{t\in D} \alpha_t^{p^{\tau-\mu}} t^{p^{\tau-\mu}}-\sum_{t\in D} \beta_t t.
\end{align}

\noindent This belongs to $M_j=M\cap k(x_1,\ldots,x_r)[x_{r+1}, \ldots, x_{l}, x_j, x_{j+1}^{p^{m}},\ldots, x_n^{p^{m}}]$, which by hypothesis equals $S_j=J\cap k(x_1,\ldots,x_r)[x_{r+1}, \ldots, x_{l}, x_j, x_{j+1}^{p^{m}},\ldots, x_n^{p^{m}}]$. By Lemma \ref{thm:L fflat mod J} we see that for each $t\in D$, there exist $\gamma_{t,t'}\in k(x_1,\ldots,x_r)$ such that $t^{p^{\tau-\mu}}=\sum_{t'\in D} \gamma_{t,t'} t'$ mod $S_j$ and $\gamma_{t,t'}$ is a ratio of polynomials having total degree at most $N_2:=\left( 2\binom{(2(n-r)(b_1(n,d))^2)^{2^n}+d + n-r}{n-r}(n-r)(b_1(n,d))^2 \right)^{2^r}$.

Substituting $\sum_{t'\in D} \gamma_{t,t'} t'$ for $t^{p^{\tau-\mu}}$, we get the following equalities mod $S_j$:

\begin{align}\label{block1}
\begin{split}
\sum_{t'\in D} \beta_{t'} t' =\sum_{t\in D} \beta_t t &=\sum_{t\in D} \alpha_t^{p^{\tau-\mu}} t^{p^{\tau-\mu}}\\
&=\sum_{t\in D} \alpha_t^{p^{\tau-\mu}}\left( \sum_{t'\in D} \gamma_{t,t'} t'\right)\\
&=\sum_{t\in D} \sum_{t'\in D}\left (\alpha_t^{p^{\tau-\mu}}  \gamma_{t,t'} t'\right)\\
&=\sum_{t'\in D} \left (\sum_{t\in D}\alpha_t^{p^{\tau-\mu}}  \gamma_{t,t'}\right) t'.
\end{split}
\end{align}

The $t_i+S_j$ form a $k(x_1,\ldots,x_r)$-basis of $A_j=B_j$, so for each $t'\in D$, we have $\sum_{t\in D}\alpha_t^{p^{\tau-\mu}}  \gamma_{t,t'}= \beta_{t'}$.

Note that the cardinality $|D|$ of $D$ is at most $N_3:=b_1(n,d)^{n-r}$ because each variable $x_{r+1}, \ldots x_j, x_{j+1}^{p^m},\ldots, x_n^{p^m}$ in $t$ has power at most $b_1(n,d)$. Clearing denominators in each equation (separately for each $t'$), we obtain a system of equations $\{\sum_{t\in D}\alpha_t^{p^{\tau-\mu}}  \widetilde{\gamma}_{t,t'}= \widetilde{\beta}_{t'}\}_{t'\in D}$ such that $\widetilde{\gamma}_{t,t'},\widetilde{\beta}_{t'}$ belong to $k[x_1,\ldots, x_r]$ and have total degree at most $N_4:= N_3N_2 +N_1$.

By Theorem \ref{thm:baby fflat bd pth power}, for each $t\in D$ there exists $\widetilde{\alpha}_t\in k(x_1,\ldots,x_r)$ that is a ratio of polynomials of total degree at most  $N_5:=(2N_3^2 (N_3+1) \binom{N_4+r}{r}p^{(\tau-\mu)(r-1)}N_4)^{2^r}$ and such that $\{\sum_{t\in D}\widetilde{\alpha}_t^{p^{\tau-\mu}}  \widetilde{\gamma}_{t,t'}= \widetilde{\beta}_{t'}\}_{t'\in D}$. Dividing to ``unclear denominators'', we have $\{\sum_{t\in D}\widetilde{\alpha}_t^{p^{\tau-\mu}}  \gamma_{t,t'}= \beta_{t'}\}_{t'\in D}$.

Referring to the equalities in (\ref{block0}) and (\ref{block1}) and substituting for $\beta_{t'}$, we obtain the following equalities mod $S_j$:

\begin{align*}\label{block1}
\begin{split}
\sum_{t\in D} \alpha_t^{p^{\tau-\mu}} t^{p^{\tau-\mu}}=\sum_{t\in D} \beta_t t\ = \sum_{t'\in D} \beta_{t'} t'&=\sum_{t'\in D} \left(\sum_{t\in D}\widetilde{\alpha}_t^{p^{\tau-\mu}}\gamma_{t,t'}\right) t'\\
&=\sum_{t'\in D} \sum_{t\in D}\left (\widetilde{\alpha}_t^{p^{\tau-\mu}} \gamma_{t,t'} t'\right)\\
&=\sum_{t\in D} \widetilde{\alpha}_t^{p^{\tau-\mu}}\left (\sum_{t'\in D}  \gamma_{t,t'}t'\right)\\
&=\sum_{t\in D} \widetilde{\alpha}_t^{p^{\tau-\mu}} t^{p^{\tau-\mu}} \\
&=\sum_{t\in D} \left(\widetilde{\alpha}_t t\right)^{p^{\tau-\mu}} .
\end{split}
\end{align*}

Since $\tau-\mu>0$ by assumption, $p^{\tau-\mu}$-th roots are unique and $\sum_{t\in D} \alpha_t t=\sum_{t\in D} \widetilde{\alpha}_t t$ mod $S_j$. It follows that $x_{j+1}^{p^{\mu}}+ \sum_{t\in D}\widetilde{\alpha}_t t\in M$. Observe that the $(x_1,\ldots,x_r)$-degree of this polynomial is at most $N_6:=N_5+(n-r)(b_1(n,d))^2$. Let $d_0$ in Corollary \ref{thm:M down to J} be $N_6$ and consequently define $\widetilde{d}$ to be $ (b_1(n,d)^{n-r})(d_0)(1+(1+b_1(n,d))^{n-r})=N_3N_6(1+(1+b_1(n,d))^{n-r})$. Since $I$ is prime up to $b_6(n,d):= max\left\{b_4(n,d), b_5(n,d), \widetilde{d}\! \left( 1+\binom{\widetilde{d}+n}{n}\right)\right\}$, by Corollary \ref{thm:M down to J} we have $x_{j+1}^{p^{\mu}}+ \sum_{t\in D}\widetilde{\alpha}_t t\in J$, in fact in $S_{j+1}=J\cap k(x_1,\ldots,x_r)[x_{r+1}, \ldots, x_{j}, x_{j+1},x_{j+2}^{p^{m}},\ldots,  x_n^{p^{m}}]$. Thus $x_{j+1}^{p^\mu}+S_{j+1}\in A_j$, contradicting the minimality of $\tau$ and proving that $S_{j+1}=M_{j+1}$. By induction, $J=S_{n}=M_{n}= M$ and $J$ is maximal.



\end{proof}


\end{subsubsection}
\end{subsection}

\begin{subsection}{Polynomial bounds on primality}
Let $b(n,d)$ be the function $b_6(n,d)$ from the preceding theorem.
\begin{theorem}\label{thm:prime}
If $I$ is prime up to $b(n,d)$, then $I$ is prime. 
\end{theorem}
\begin{proof}
The cases $n=1$ and $d=1$ are covered by \ref{thm:easy cases}. If $n,d>1$, then $J\cap k[x_1,\ldots,x_n] =I$ by \ref{thm:J intersect}. The ideal $J$ is maximal by \ref{thm:maximal} and restrictions of prime ideals to subrings are prime, so $I$ is prime.
\end{proof}

We conclude this section with concrete upper bounds on $b(n,d)$. For each lemma the claimed bound is greater than $d$, so the cases $n=1$ and $d=1$ are covered and we may assume $n,d>1$ where necessary.

\begin{lemma}
$b_4(n,d)\leq (2nd)^{n2^{4n+8}}$ for all $n,d\geq 1$.
\end{lemma}
\begin{proof}
Using earlier bounds on $b_1(n,d)$ and $b_3(n,d)$ we obtain:
\begin{align*}
&b_4(n,d)=b_1(n+1, (n-1)\cdot b_3(n,d)+1)\leq  2((n-1)\cdot b_3(n,d)+1)^{2^{n+1}} \\ &\leq 2((n-1)\cdot 2^{n2^{2n+6}}d^{n2^{3n+2}}
+1)^{2^{n+1}} \leq 2(n2^{n2^{2n+6}}d^{n2^{3n+2}}
)^{2^{n+1}}  \\&
\leq (n2^{n2^{2n+7}}d^{n2^{3n+2}}
)^{2^{n+1}} \leq  (2nd)^{n2^{4n+8}}. 
\end{align*}
\end{proof}


Similarly we obtain bounds on the values defined in \ref{def:intermed bds}:
\begin{lemma}
$N_1\leq N_2\leq (2nd)^{n2^{3n+2}}$.
\end{lemma}
\begin{proof}
Clearly $N_1\leq N_2$.  For the second inequality, compute

\begin{align*}
&N_2 = \left( 2\binom{(2(n-r)(b_1(n,d))^2)^{2^n}+d + n-r}{n-r}(n-r)(b_1(n,d))^2 \right)^{2^r} \\& \leq
 \left( 2\binom{(2n(b_1(n,d))^2)^{2^n}+d + n}{n}(n)(b_1(n,d))^2 \right)^{2^n} \\&\leq
\left( 2\binom{(2n(2d^{2^n})^2)^{2^n}+d + n}{n}(n)(2d^{2^n})^2 \right)^{2^n}\\&\leq
\left( 2\left((2n(2d^{2^n})^2)^{2^n}+d + n\right)^n(n)(2d^{2^n})^2 \right)^{2^n} \\& \leq
\left( 2\left(2(2n(2d^{2^n})^2)^{2^n}\right)^n(n)(2d^{2^n})^2 \right)^{2^n}\\&\leq
2^{5n2^{2n+1}}n^{2n2^{2n}}d^{2n2^{3n+1}}\\& \leq
(2nd)^{n2^{3n+2}}.
\end{align*}
\end{proof}

\begin{lemma}
$N_3\leq 2^nd^{n2^{n}}$.
\end{lemma}
\begin{proof}
\begin{align*}
&N_3= b_1(n,d)^{n-r}\leq b_1(n,d)^n \leq (2d^{2^n})^n =  2^nd^{n2^{n}}.
\end{align*}
\end{proof}

\begin{lemma}
$N_4\leq (2nd)^{n2^{3n+3}}$.
\end{lemma}
\begin{proof}
\begin{align*}
&N_4= N_3N_2 +N_1\leq N_2(N_3+1)\leq (2nd)^{n2^{3n+2}}( 2^nd^{n2^{n}}+1)\\&
\leq (2nd)^{n2^{3n+2}}( 2^{n+1}d^{n2^{n}})\leq (2nd)^{n2^{3n+3}}.
\end{align*}
\end{proof}

\begin{lemma}
$N_5\leq (2nd)^{n(n+1)2^{4n+4}}$.
\end{lemma}
\begin{proof}

\begin{align*}
&N_5 = (2N_3^2 (N_3+1) \binom{N_4+r}{r}(b_1(n,d))^{(r-1)}N_4)^{2^r} \\& 
\leq \left(2^2N_3^3  \binom{N_4+n}{n}(2d^{2^n})^{(n-1)}N_4)\right)^{2^n} \\& \leq
\left(2^2N_3^3 (N_4+n)^n(2d^{2^n})^{(n-1)}N_4)\right)^{2^n} \\&
\leq \left(2^2N_3^3 (2N_4)^n(2d^{2^n})^{(n-1)}N_4)\right)^{2^n}\\&
\leq \left(2^2(2^nd^{n2^{n}})^3 (2( (2nd)^{n2^{3n+3}}))^n(2d^{2^n})^{(n-1)}(2nd)^{n2^{3n+3}}\right)^{2^n}\\&
\leq 2^{(1+2n)2^n}(2^nd^{n2^{n}})^{3\cdot 2^n} ((2nd)^{n2^{3n+3}})^{(n+1)2^n}d^{(n-1)2^{2n}}\\&
\leq 2^{(1+2n)2^n+ 3n2^n +n(n+1)2^{4n+3}}n^{n(n+1)2^{4n+3}}d^{3n2^{2n}+n(n+1)2^{4n+3}+(n-1)2^{2n}}\\&
\leq  (2nd)^{(4n-1)2^{2n} + n(n+1)2^{4n+3}} \leq  (2nd)^{n(n+1)2^{4n+4}}.
\end{align*}
\end{proof}

\begin{lemma}
$N_6\leq 2(2nd)^{n(n+1)2^{4n+4}}$.
\end{lemma}
\begin{proof}
\begin{align*}
&N_6= N_5+(n-r)(b_1(n,d))^2\leq N_5+n(2d^{2^n})^2\\& \leq  (2nd)^{n(n+1)2^{4n+4}}+ n2^2d^{2^{n+1}}\leq 2(2nd)^{n(n+1)2^{4n+4}}.
\end{align*}
\end{proof}

\begin{lemma}
$\widetilde{d}\leq (2nd)^{n(n+1)2^{4n+5}}$.
\end{lemma}
\begin{proof}
\begin{align*}
&\widetilde{d}=N_3N_6(1+(1+b_1(n,d))^{n-r})\leq N_3N_6(1+(1+2d^{2^n})^n)\leq  N_3N_6(2^{2n+1}d^{n2^n})\\&
\leq (2^nd^{n2^{n}})(2(2nd)^{n(n+1)2^{4n+4}})(2^{2n+1}d^{n2^n})\\& 
\leq 2^{3n+2 +n(n+1)2^{4n+4} }n^{n(n+1)2^{4n+4}}d^{2n2^n +n(n+1)2^{4n+4}}\\&
\leq (2nd)^{n(n+1)2^{4n+5}}.
\end{align*}
\end{proof}

\begin{lemma}
$\widetilde{d} \left( 1+\binom{\widetilde{d}+n}{n}\right)\leq (2nd)^{n(n+1)^22^{4n+6}} $.
\end{lemma}
\begin{proof}
\begin{align*}
&\widetilde{d} \left( 1+\binom{\widetilde{d}+n}{n}\right) \leq \widetilde{d}\left(1+(\widetilde{d}+n)^n\right) \leq \widetilde{d}(1+(2\widetilde{d})^n)\leq (2\widetilde{d})^{n+1}\\& 
\leq \left(2(2nd)^{n(n+1)2^{4n+5}} \right) ^{n+1} \leq (2nd)^{n(n+1)^22^{4n+6}}.
\end{align*}
\end{proof}

\begin{theorem}\label{thm:final estimate}
$b(n,d)\leq max\left\{ (2d)^{2^{3n^2+2n}}, (2nd)^{n(n+1)^2{2^{4n+6}}}\right\} $ for all $n,d\geq~\!1$.
\end{theorem}
\begin{proof}

\begin{align*}
&b(n,d)=b_6(n,d)=max\left\{b_4(n,d), b_5(n,d), \widetilde{d}\! \left( 1+\binom{\widetilde{d}+n}{n}\right)\right\}\\&
\leq max \left\{ (2nd)^{n2^{4n+8}},(2d)^{2^{3n^2+2n}},  (2nd)^{n(n+1)^22^{4n+6}} \right\}\\&
\leq  max\left\{ (2d)^{2^{3n^2+2n}}, (2nd)^{n(n+1)^2{2^{4n+6}}}\right\}.
\end{align*}
\end{proof}

If desired, one can continue to simplify and obtain larger, more-readable bounds: 

\begin{cor}\label{thm:simple bd}
$b(n,d) \leq (2nd)^{n^32^{6n^2}}$ for all $n,d~\geq~1$.
\end{cor}

\begin{proof}

\begin{align*}
&b(n,d)\leq max\left\{ (2d)^{2^{3n^2+2n}}, (2nd)^{n(n+1)^2{2^{4n+6}}}\right\}\\& 
\leq (2nd)^{max\left\{ 2^{3n^2+2n}, n(n+1)^2{2^{4n+6}}\right\} }\leq (2nd)^{n(n+1)^22^{3n^2+2n}}\\&
\leq (2nd)^{n^32^{3n^2 +2n+2}}\leq (2nd)^{n^32^{6n^2}}.
\end{align*}
\end{proof}
\end{subsection}
\end{section}

\section{Maximal Ideals}

We note that an analogous result, with a simpler proof, holds for maximal ideals: if an ideal is \emph{maximal up to} $b$ (see \ref{def:bdd prime and max}) for large enough $b$, then the ideal is maximal. Similar remarks apply to the following proof and an ultraproduct argument of Schoutens (4.1.4, \cite{schoutens2010use}) as held for Theorem \ref{thm:prime bound} and Schmidt-G\"ottsch's argument mentioned in the introduction.

The following lemma is a standard step in the proof of the Noether Normalization Theorem (see, e.g., 1.1.8 in \cite{crespo2011algebraic}):
\begin{lemma}\label{thm:noether lemma}
  Let $f\in k[x_1,\ldots,x_n]$ be a nonzero polynomial of total degree $d$ over a field $k$.
There is some $c\in k\setminus{\{0\}}$ and $a_1,\ldots,a_{n-1}\leq (d+1)^{n-1}$ such that the polynomial
  \[cf(y_1+y_n^{a_1},y_2+y_n^{a_2},\ldots,y_{n-1}+y_n^{a_{n-1}},y_n)\]
  is monic in $y_n$.
\end{lemma}
To streamline the main proof we first take care of the simple cases and prove a technical inequality.
\begin{lemma}\label{thm:lin max bd}
For any field $k$ and any proper ideal $I\subseteq k[x_1,\ldots,x_n]$ with generators of degree $1$, if $I$ is maximal up to $1$, then $I$ is maximal. If $n=1$ and $I$ is proper with generators of degree at most $d\geq 1$, then maximality up to $d$ suffices. 
\end{lemma}
\begin{proof}
Suppose the generators of $I$ have degree 1. By an invertible affine change of variables, we may assume that the generators have the form $x_1,\ldots, x_m$ for some $m\leq n$. (Such a change of variables preserves our assumption that linear polynomials not in the ideal are invertible modulo the ideal.) If $m<n$, then $x_{m+1}$ is invertible mod I and $x_{m+1}f -1\in I$ for some $f$. However, $I=(x_1,\ldots ,x_m)$ cannot have an element with nonzero constant term. Hence $m=n$ and $I$ is maximal.

In the case $n=1$, $k[x_1]$ is a PID and $I=(f)$ for some $f$ of degree at most $d$. If $f=0$, then $x$ is not invertible mod $I$ and $I$ is not maximal up to $d$. Otherwise, reduction by $f$ gives a nonzero remainder of degree less than $d$ for every element of $k[x_1]\setminus I$. Hence maximality up to $d$ implies that $I$ is maximal.
\end{proof}

\begin{lemma}\label{thm:max bd ind}
Let $k,n,d$ be natural numbers such that $3<k+1<n$ and $d\geq 1$. Define $D$ to be $(2d)^{n(n-1)\cdots(n-k)2^{(k-1)n}}$. 
Then $(2D)^{(n-k-1)2^{(n-k+2)}}\leq (2d)^{n (n-1)\cdots (n-k)(n-k-1) 2^{kn}}$.
\end{lemma}

\begin{proof}
Applying the definition of $D$ and simplifying, we convert the left side of the claimed inequality into 

\[2^{(n-k-1)2^{n-k+2}}(2d)^{n\cdots(n-k-1)2^{kn-k+2}}.\] 

Dividing by $(2d)^{n\cdots(n-k-1)2^{kn-k+2}}$ and simplifying, we see the claimed inequality is equivalent to  

\begin{align*}
2^{(n-k-1)2^{n-k+2}}&\leq \left((2d)^{n (n-1)\cdots (n-k)(n-k-1)}\right)^{( 2^{kn}-2^{kn-k+2})}\\
&= \left((2d)^{n (n-1)\cdots (n-k)(n-k-1)}\right)^{( 2^{kn-k+2})(2^{k-2}-1)},
\end{align*}

\noindent which is clear since $k>2$.
\end{proof}

\begin{theorem}\label{thm:bdd max}
Define $m(n,d):=(2d)^{(n-1)(n!)2^{(n-2)n}}$ for all $n>2,d> 1$. Otherwise let $m(n,1):=1$ (for any $n$) and $m(1,d):=d, m(2,d):=4d^4$ if $d>1$.  For any field $k$ and any proper ideal $I\subseteq k[x_1,\ldots,x_n]$ with generators of total degree at most $d$, if $I$ is maximal up to $m(n,d)$, then $I$ is maximal.  



\end{theorem}
\begin{proof}


By Lemma \ref{thm:lin max bd}, bounds of $m(n,1)=1$ and $m(1,d)=d$  suffice if $d=1$ or $n=1$ and $d>1$, respectively. Thus we may assume $n,d>1$; by Theorem \ref{thm:GB Dube} $b_1(n,d)\leq 2d^{2^n}$. 

The Noether Normalization Theorem applied to $R=k[x_1,\ldots,x_n]/I$ states that there is an injective homomorphism $\pi:k[z_1,\ldots,z_r]\rightarrow R$ such that $R$ is a finite extension of a polynomial ring $k[z_1,\ldots,z_r]$ for some $0\leq r\leq n$. The usual proof proceeds by induction on the number of generators $n$ of $R$ as a $k$-algebra. We reproduce the argument but keep track of bounds along the way. Our first goal is to show that in our scenario, bounded maximality forces $r$ to be 0. 

Toward a contradiction, suppose $r>0$. By maximality up to $m(n,d)$, the ideal $I$ must be nonzero.  Starting with a nonzero generator of degree at most $d$, Lemma \ref{thm:noether lemma} changes variables by substituting $y_i+ y_n^{a_i}$ for $x_i$ if $i<n$ and $y_n$ for $x_n$ using some integers $a_i\leq (d+1)^{n-1}$ such that the transformed generator is monic in $y_n$. This yields generators of degree at most $(d+1)^{n-1} \cdot d \leq (d+1)^n$ for the ideal $\widetilde{I}$ in $k[y_1,\ldots, y_{n}]$ generated by the generators of $I$ (after substitution). The induction continues if there is a nontrivial relation between $y_1,\ldots, y_{n-1}$; i.e., if the elimination ideal  obtained by eliminating $y_n$ from $\widetilde{I}$ is nonzero. Otherwise $\pi(z_i)= (x_i-x_n^{a_i})+ I$ for $i<n$ gives the claimed map; note that $x_n+I$ is integral over $\pi(k[z_1,\ldots,x_{n-1}])$.  By hypothesis, the process stops after at most $n-1$ steps, leaving us with the desired map for some $r\geq 1$ and $k[z_1,\ldots,z_r]$. Theorem \ref{thm:GB elim bd} implies that the generators of the elimination ideal are bounded by $b_1(n,(d+1)^{n})$ after one step of substitution and elimination, $b_1(n-1,(b_1(n,(d+1)^{n})+1)^{n-1})$ after two, and so on. 

 We analyze the alternating substitution and elimination steps to bound the degrees of (representatives of) $\pi(z_1),\ldots, \pi(z_r)$ as polynomials in the variables $x_1,\ldots,x_n$. Let $E_n=d$. Define $D_{n-k}$ to be $(E_{n-k+1}+1)^{n-k}$ and $E_{n-k}$ to be $2(D_{n-k}\cdot E_{n-k+1})^{2^{n-k+1}}$ for $0<k<n$. By induction, $D_{n-k}$ is a bound on the degree of a substitution at the $k$-th step and $E_{n-k}$ is a bound on the degree of the generators after $k$-many steps of substitution and elimination. Since we seek the degrees of the images of $z_1,\ldots,z_r$ in the original variables (and not just the variables in the preceding step), our goal is to bound the product $\prod_{i=1}^{n-1} D_{n-i}$. (We eliminate at most $n-1$ variables lest $r=0$.) 

Note that $E_{n-k+1}\leq D_{n-k}$ and so $E_{n-k}\leq 2(D_{n-k}^2)^{2^{n-k+1}}= 2(D_{n-k})^{2^{n-k+2}}$ for all $0<k<n.$ It follows that $E_{n-k+1}\leq 2(D_{n-k+1})^{2^{n-k+3}}$ and hence 

\[D_{n-k}\leq (2(D_{n-k+1})^{2^{n-k+3}}+1)^{n-k}\leq (2^2(D_{n-k+1})^{2^{n-k+3}})^{n-k}\leq (2D_{n-k+1})^{(n-k)2^{n-k+3}}\]

\noindent for $1<k<n$. 

It is immediate from the definitions that $D_{n-1}=(d+1)^{n-1}\leq (2d)^{n(n-1)}$, and a calculation like that of Lemma \ref{thm:max bd ind} shows that $D_{n-2}\leq (2d)^{n(n-1)(n-2)2^{n}}$. Using Lemma \ref{thm:max bd ind} and the inequality $D_{n-k}\leq (2D_{n-k+1})^{(n-k)2^{n-k+3}}$ for induction on $k$, we find that $D_{n-k}\leq (2d)^{n(n-1)\cdots (n-k)2^{(k-1)n}}$ for $2<k<n$; the previous sentence establishes the inequality for $k=1,2$. It follows that $\prod_{i=1}^{k} D_{n-i} \leq (2d)^{kn(n-1)\cdots(n-k)2^{(k-1)n}}$ and so $\pi(z_1),\ldots,\pi(z_r)$ can be represented by polynomials of degree at most $(2d)^{(n-1)(n!)2^{(n-2)n}}$ in $x_1,\ldots,x_n$. One easily checks that $(2d)^{(n-1)(n!)2^{(n-2)n}}\leq m(n,d)$ for all $n,d> 1$. 


Then since $I$ is maximal up to $m(n,d)$, $\pi(z_r)$ is invertible in $R$. ($z_r\neq 0$, so by injectivity of $\pi$ we know $\pi(z_r)\notin I$.) Since $R$ is integral over $\pi(k[z_1,\ldots, z_r])$, $(\pi(z_r))^{-1}$ is a root of a monic equation with coefficients in $\pi(k[z_1,\ldots, z_r])$. Multiplying the monic equation by an appropriate power of $\pi(z_r)$ shows that $(\pi(z_r))^{-1}\in  \pi(k[z_1,\ldots, z_r])$, whence $z_r$ is invertible in $k[z_1,\ldots,z_r]$. But this is absurd, so $r$ is equal to 0 and $R$ is finite over $k$.

Since $R$ is finite over $k$, for each $x_1,\ldots, x_n$ the restricted ideal $I\cap k[x_i]$ is nonzero. 
By 
Theorem \ref{thm:GB elim bd} there is a nonzero polynomial $\alpha_i(x_i)\in I\cap k[x_i]$ of degree at most $b_1(n,d)$ for each $x_1,\ldots, x_n$. Using $\alpha_i$ to reduce higher powers of $x_i$ mod $I$, we see that any $f\notin I$ is equivalent mod $I$ to some $\widetilde{f}\notin I$ of total degree at most $nb_1(n,d)\leq m(n,d)$ for $n,d>1$.
 By maximality up to $m(n,d)$, $f$ is invertible mod $I$  and hence $I$ is maximal. 


\end{proof}
A single formula for the bound easily follows: 

\begin{cor}
$m(n,d)\leq (2d)^{(n)(n!)2^{(n-1)n}}$ for all $n,d\geq 1$.
\end{cor}

\bibliographystyle{plain}
\bibliography{PolyBounds} 
\end{document}